\documentclass[11pt,amsfonts]{article}
\usepackage{graphicx}
\usepackage{latexsym}
\usepackage{amssymb}
\usepackage{amsmath}
\usepackage{enumerate}
\usepackage{color}
\usepackage{layout}
\usepackage{dsfont}
\usepackage{eufrak}
\newtheorem{prop}{Proposition}
\newtheorem{lemma}{Lemma}
\newtheorem{definition}{Definition}
\newtheorem{corollary}{Corollary}
\newtheorem{theorem}{Theorem}
\newtheorem{remark}{Remark}

\def\real{{\mathord{{\rm I\kern-2.8pt R}}}}        
\def\inte{{\mathord{{\rm I\kern-2.8pt N}}}}

\def\sZZ{{\rm Z\kern-2.8ptem{}Z}}

\def\z{{\mathchoice
  {\sZZ}
  {\sZZ}
  {\rm Z\kern-0.30em{}Z}
  {\rm Z\kern-0.25em{}Z} }}
\def\sQQ{{\kern 0.27em \vrule height1.45ex width0.03em depth0em
          \kern-0.30em \rm Q}}
\def\qu{{\mathchoice
    {\sQQ}
    {\sQQ}
  {\kern 0.225em \vrule height1.05ex width0.025em depth0em \kern-0.25em \rm Q}
  {\kern 0.180em \vrule height0.78ex width0.020em depth0em \kern-0.20em \rm Q}
        }}
\def\sCC{{\kern 0.27em \vrule height1.45ex width0.03em depth0em
          \kern-0.30em \rm C}}
\def\complex{{\mathchoice
    {\sCC}
    {\sCC}
  {\kern 0.225em \vrule height1.05ex width0.025em depth0em \kern-0.25em \rm C}
  {\kern 0.180em \vrule height0.78ex width0.020em depth0em \kern-0.20em \rm C}
        }}

\newcommand{\ba}{\begin{array}}
\newcommand{\ea}{\end{array}}
\newcommand{\be}{\begin{equation}}
\newcommand{\ee}{\end{equation}}
\newcommand{\bea}{\begin{eqnarray}}
\newcommand{\eea}{\end{eqnarray}}
\newcommand{\beaa}{\begin{eqnarray*}}
\newcommand{\eeaa}{\end{eqnarray*}}

\def\z{\zeta}

\font\tenmath=msbm10 \font\sevenmath=msbm7 \font\fivemath=msbm5
\newfam\mathfam \textfont\mathfam=\tenmath
\scriptfont\mathfam=\sevenmath \scriptscriptfont\mathfam=\fivemath

\def \={{\buildrel {\rm (law)} \over =}}

\def\L{\Lambda}

\def\qed{ \hfill \vrule width.25cm height.25cm depth0cm\smallskip}

\newcommand{\basa}{\begin{assumption}}
\newcommand{\easa}{\end{assumption}}

\newcommand{\bas}{\begin{assum}}
\newcommand{\eas}{\end{assum}}

\def\limsup{\mathop{\overline{\rm lim}}}

\def\bea{\begin{eqnarray}}
\def\eea{\end{eqnarray}}
\def\beas{\begin{eqnarray*}}
\def\eeas{\end{eqnarray*}}
\def\bd{\begin{description}}
\def\ed{\end{description}}
\def\im{\item}

\def\bd{\begin{description}}
\def\ed{\end{description}}

\def\ba{\bar{A}}

\def\L{{\bf L}}

\def\be{\begin{equation}}
\def\ee{\end{equation}}
\def\bea{\begin{eqnarray}}
\def\eea{\end{eqnarray}}
\def\beas{\begin{eqnarray*}}
\def\eeas{\end{eqnarray*}}
\def\bi{\begin{itemize}}
\def\ei{\end{itemize}}
\def\im{\item}
\def\bd{\begin{description}}
\def\ed{\end{description}}
\def\r{\right}

\def\L{\Lambda}

\newcommand{\bbD}{{\mathbb D}}

\newcommand{\bbN}{{\mathbb N}}

\newcommand{\bbR}{{\mathbb R}}

\newcommand{\bbZ}{{\mathbb Z}}

\newcommand{\ignore}[1]{}
\begin{document}

\renewcommand{\thefootnote}{\fnsymbol{footnote}}

\renewcommand{\thefootnote}{\fnsymbol{footnote}}

\title{Asymptotic expansion for vector-valued sequences of random variables with focus on Wiener chaos }
\author{ Ciprian A. Tudor $^{1 }$ \footnote{. The author was partially supported by Project ECOS  NCRS-CONICYT C15E05,MEC PAI80160046 and MATHAMSUD 16-MATH-03 SIDRE Project.} \qquad  Nakahiro Yoshida  $^{2,3,}$\footnote{The author was in part supported by Japan Science and Technology Agency CREST JPMJCR14D7; Japan Society for the Promotion of Science Grants-in-Aid for Scientific Research No. 17H01702 (Scientific Research) and by a Cooperative Research Program of the Institute of Statistical Mathematics.  } \vspace*{0.1in} \\
$^{1}$ Laboratoire Paul Painlev\'e, Universit\'e de Lille 1\\
 F-59655 Villeneuve d'Ascq, France.\\
\quad tudor@math.univ-lille1.fr\vspace*{0.1in}\\
$^{2}$Graduate School of Mathematical Science, University of Tokyo\\
 3-8-1 Komaba, Meguro-ku, Tokyo  153, Japan \\
 $^{3}$ Japan Science and Technology Agency CREST\\
 \quad nakahiro@ms.utokyo.ac.jp \vspace*{0.1in} }
\maketitle

\begin{abstract}
We develop the asymptotic expansion theory for vector-valued sequences $(F_{N})_{N\geq 1}$ of random variables in terms of the convergence of the Stein-Malliavin matrix  associated to the sequence $F_{N}$. Our approach combines the classical Fourier approach and the recent theory on Stein method and Malliavin calculus. We find the second order term of the asymptotic expansion of the density of $F_{N}$ and we illustrate our results by several examples.
\end{abstract}

\vskip0.3cm

{\bf 2010 AMS Classification Numbers: } 62M09, 60F05, 62H12

\vskip0.3cm


\section{Introduction}

The analysis of the convergence in distribution of sequences of random variables constitutes a fundamental  direction in probability and statistics. In particular, the asympotic expansion represents a technique that allows to give a precise approximation for the densities of probability distributions,
 
Our work concerns the convergence in law of random sequences  to the Gaussian distribution. Although our context is more general, we will focus on examples based on elements in Wiener chaos.  The characterization of the convergence in distribution  of sequences of random variables belonging to a Wiener chaos of fixed order  to the Gaussian law constitutes an important research direction in stochastic analysis in the last years. This research line was initiated by the seminal paper \cite{NuPe}  where the authors proved the famous Fourth Moment Theorem. This result, which basically says that the convergence in law of a sequence of multiple stochastic integrals with unit variance  is equivalent   to  the convergence of the sequence of the fourth moments to the fourth moment of the Gaussian distribution,  has then been extended, refined, completed or applied to various situations. The reader may consult the recent monograph \cite{NPbook} for the basics of this theory. 

In this work, we are also concerned with some refinements of the Fourt Moment Theorem for sequences of random variables in a Wiener chaos of fixed order. More concretely, we aim to find  the second order  expansion for   the density.  This is a natural extension of the Fourth Moment Theorem and of its ramifications. In order to explain this link, let us briefly describe the context. 
Let $H$ be a real and separable Hilbert space, $(W(h), h \in H)$ an isonormal Gaussian process on a probability space $\left( \Omega, \mathcal{F}, P\right)$ and let $ I_{q}$ denote the $q$th multiple stochastic  integral with respect to $W$. Let $Z\sim N(0,1)$. Denote by $d_{TV}(F,G)$ the total variation distance between  the laws of the random variables $F$ and $G$ and let $D$ be the Malliavin derivative with respect to $W$.  The Fourth Moment Theorem from  \cite{NuPe} can be stated as follows.
\begin{theorem}
Fix $q \geq 1$. Consider a sequence $\{ F_k=I_q(f_k), k\geq 1\}$ of  random variables in the $q$-th Wiener chaos.  Assume that
\begin{equation}
\lim _{k\rightarrow \infty} \mathbf{E}F_k^2 = \lim _{k\rightarrow \infty}q! \| f_k\| ^2_{H^{\otimes q}} =1.
\end{equation}
Then, the following statements are equivalent:
\begin{enumerate}
\item \label{4momO1} The sequence of random variables $\{ F_k=I_n(f_k), k\geq 1\}$ converges  in distribution as $k\rightarrow \infty$ to the standard normal law $N (0,1)$.
\item \label{4momO2} $\lim _{k\rightarrow \infty}\mathbf{E}[F_k^4]=3$.

\item \label{4momO4} $\| DF_k\| _H^2$ converges to $q$ in $L^2 (\Omega)$ as $k\rightarrow \infty$.

\item $d_{TV} (F_{k}, Z)\to _{k\to \infty} 0.$\label{4}
\end{enumerate}
\end{theorem}

A more recent point of interest in the analysis of  the convergence to the normal distribution of sequences of multiple integrals is represented by the analysis of the convergence of the sequence of densities. For every $N\geq 1$, denote by $p_{F_{N}}$ the density of the random variable $F_{N}$ (which exists due to a result in \cite{Shige}). Since the total variation distance between $F$ and $G$ (or between the distributions of the random variables $F$ and $G$) can be written as $d_{TV} (F,G)= \frac{1}{2} \int_{\mathbb{R} }   \vert p_{F}(x) -p_{G}(x) \vert dx, $ we can notice that point \ref{4}. in the Fourth Moment Theorem implies that 
\begin{equation*}
\Vert p_{F_{N}} (x)- p(x) \Vert _{L^{1} (\mathbb{R})} \to _{N\to \infty} 0
\end{equation*}
where $p(x)= \frac{1}{\sqrt{2\pi}} e ^{-\frac{x ^{2}}{2}}$ denotes, throughout our work, the density of the standard normal distribution. The result has been improved  in the work \cite{NoNu}, where  the authors showed that if 
\begin{equation}
\label{neg}
\limsup _{N\to \infty} \mathbf{E}\Vert DF_{N} \Vert _{H}^{-4-\varepsilon} <\infty
\end{equation}
then the conditions \ref{4momO1}.-\ref{4}. in the Fourth Moment Theorem  are equivalent to

\begin{equation}\label{24i-1}
\sup_{x\in \mathbb{R}} \left| p_{F_{N}} (x)- p(x) \right| \leq c\sqrt{ (\mathbf{E}F_{N}^{4}- 3)}.
\end{equation}
We also refer to \cite{HLN}, \cite{HNTX} for other references related to the density convergence on Wiener chaos.

Our purpose is to find the asymptotic expansion in the Fourth Moment Theorem (and, more generally, for sequences of random variables converging to the normal law). By finding the asymptotic expansion we mean  to find a sequence $(p_{N}) _{N\geq 1}$ of the form  $p_{N}(x)=p(x)+ r_{N} \varphi (x) $ such that  $r_{N}\to _{N\to \infty} 0$ which approximates the density $p_{F_{N}}$ in the following sense

 \begin{equation}
\label{12i-1}\sup_{x\in \mathbb{R}} \vert x ^{\alpha}\vert \left| p_{F_{N}}(x)- p_{N}(x) \right| = o(r_{N})
\end{equation}
for any $\alpha \in \mathbb{Z}_{+}.$
This improves the main result in \cite{NoNu} and it will also allow to generalize the Edgeworth -type expansion (see \cite{NoPe2}, \cite{BBNP}) for sequences on Wiener chaos to a larger class of measurable functions. The main idea to get the asymptotic expansion is to analyze the asymptotic behavior of the characteristic function of $F_{N}$ and then to obtain the approximate density $p_{N}$ by Fourier inversion of the dominant part of the characteristic function.

As mentioned before, we develop our asymptotic expansion theory for a general class of random variables, which are not necessarily multiple stochastic integrals, although our main examples are related to Wiener chaos.  Our main assumption is that the vector-valued  random sequence $\left(F_{N}, \gamma_{N}^{-1} (\Lambda (F_{N})-C) \right)$ converges to a $d+d\times d$-dimensional Gaussian vector, where $(\gamma_{N})_{N\geq 1}$ is a deterministic sequence converging to zero as $N\to \infty$,  $\Lambda (F_{N})$ is the so-called Stein-Malliavin matrix of $F_{N}$ whose components are $\Lambda _{i, j} =\langle DF_{N} ^{(i) }, D(-L) ^{-1} F_{N} ^{(j)}\rangle _{H}$  for $i,j=1,d$ ($L$ is the Ornstein-Uhlenbeck operator)  and $C_{i,j}=\lim_{N\to \infty} \mathbf{E} F_{N} ^{(i)}F_{N} ^{(j)}$.  In the one-dimensional case the assumption becomes 
\begin{equation}
\label{intro11}
\left(F_{N}, \gamma _{N}^{-1} (\langle DF_{N}, D(-L) ^{-1}F_{N}\rangle _{H}-1)\right)
\end{equation} converges in distribution, as $N\to \infty$, to a two-dimensional centered Gaussian vector $(Z_{1}, Z_{2})$ with $ \mathbf{E }Z_{1}^{2}=\mathbf{ E} Z_{2} ^{2}=1$ and $\mathbf{E}Z_{1}Z_{2}=\rho \in (-1, 1)$.

The condition (\ref{intro11}) is related to what is usually assumed in the asymptotic expansion theory for martingales, which has  been developed in the last decades and it is based on the so-called Fourier approach.  Let us refer, among many others, to \cite{My}, \cite{Yos1}, \cite{PY}, \cite{Yos2}, \cite{Yos3} for few important works on the martingale approach in asymptotic expansion. Recall that if $(M_{N}) _{N\geq 1}$ denotes a  sequence of random variables converging to $Z\sim N(0,1)$ such that  $M_{N}$ is the terminal value of a continuous martingale and if $\langle M \rangle  _{N} $ is the martingale's bracket,  then one assumes in e.g.  \cite{Yos1}), \cite{Yos2} that  

$$\langle M\rangle _{N} \to _{N\to \infty} 1  \mbox{ in probability} $$
 and the following joint convergence holds true
\begin{equation}\label{12i-2}
\left( M_{N}, \frac{ \langle M \rangle _{N} -1}{ \gamma_{N}}\right) \to _{N\to \infty}  ^{(d)}  (Z_{1}, Z_{2}) 
\end{equation}
where $(Z_{1}, Z_{2}) $ is a centered Gaussian vector as above. The notation "$\to ^{(d) } $" stands for the convergence in distribution. 

In our assumption (\ref{intro11}), the role of the martingale bracket is played by \\ $\langle DF_{N}, D(-L) ^{-1}F_{N}\rangle _{H}$ which becomes $\frac{1}{q} \Vert DF_{N}\Vert _{H} ^{2}$ when $F_{N}$ is an element of the  $q$th Wiener chaos. The choice is natural, suggested by the identity $\mathbf{E}F_{N}  ^{2}= \mathbf{E}\langle F \rangle _{N} ^{2}$ which playes the role of the It\^o isometry and by the fact that $\frac{\Vert DF_{N}\Vert ^{2}}{q}-1$  converges to zero in $ L^{2}(\Omega)$, due to the Fourth Moment Theorem. 

We organized our paper as follows.  In Section 2, we fix the general context of our work, the main asumptions and some notation. In Section 3, we analyze the asymptotic behavior of the characteristic function of a random sequence that satisfies our assumptions while in Section 4 we obtain the approximate demsity and the asymptotic behavior via Fourier inversion. In Section 5, we illustrate our theory by several examples related to Wiener chaos.

\section{Asssumptions and notation}
We will here present the basic assumptions and the  notation  utilised throughout our work.

Let $H$ be  a real and separable Hilbert space endowed with the inner product $\langle \cdot, \cdot \rangle_{H}$ and consider $(W(h), h \in H)$ an isonormal process on the probability space $(\Omega, \mathcal{F}, P)$. In the sequel,  $D, L$ and $\delta $ represent the Malliavin derivative, the Ornstein-Uhlenbeck operator and the divergence integral with respect to $W$. See the Appendix for their definition and basic properties.

\subsection{The main assumptions}\label{21}
Consider a sequence of centered  random variables $(F_{N}) _{N\geq 1}$ in $\mathbb{R} ^{d}$ of the form
$$F_{N}= \left( F_{N} ^{(1)}, \ldots, F_{N} ^{(d)}\right).$$
  Denote by $\Lambda (F_{N})  = (\Lambda_{i,j}) _{i,j=1,..,d}$ the Stein-Malliavin matrix with components
\begin{equation*}
\Lambda _{i,j}= \langle DF_{N} ^{(i) }, D(-L) ^{-1} F_{N} ^{(j)}\rangle _{H}
\end{equation*}
for every $i,j=1,..,d$.

We consider the following assumptions:

\bd
\im[[C1\!\!]] There exists a symmetric (strictly) positive definite $d\times d$ matrix $C$ and a deterministic sequence $(\gamma_{N}) _{N\geq 1}$ converging to zero as $N\to \infty$ such that 
$\big(F_N,\gamma _N^{-1}  (\Lambda (F_{N}) -C    )\big)$ converges in distribution to 
a  $d+d\times d$  Gaussian random vector $(Z_1,Z_2)$  such that $Z_{1}= (Z_{1} ^{(1)}, \ldots , Z_{1} ^{(d)})$, $Z_{2} = (Z_{2} ^{(k,l) })_{k,l=1,.., d} $ with
\begin{equation*}
\mathbf{E} Z_{1} ^{(i)} Z_{1}^{(j)} =C_{i,j}, \mbox{ for every } i,j=1,.., d.
\end{equation*}
 \ed
\bd
\im[[C2\!\!]] For every $i=1,..,d$, the sequence 
$(F_N ^{(i)})_{N\geq 1}$ is bounded in $\bbD ^ {\ell+1,\infty}=\cap_{p>1}\bbD ^{\ell+1,p}$ for some $l\geq d+3$. 
\ed
\bd
\im[[C3\!\!]] With $\gamma_{N}, C$ from {\bf  [C1\!\! ]},  for every $i,j=1,d$ the sequence
$\left(\gamma_{N}^{-1}\big(\Lambda _{i,j}-C_{i,j} \big)\right)_{N\in\bbN}$ is bounded in 
$\bbD ^ {\ell,\infty}$ for some $l\geq d+3$. 
\ed

Notice that the matrix $\Lambda (F_{N})$ is not necessarily symmetric (except when all the components belong to the same Wiener chaos) since in general  $\Lambda _{i,j}\not=\Lambda _{j,i} $. Nevertheless, we assumed in {\bf [C1\!\! ]} that it converges to the symmetric matrix $C$, in the sense that $\gamma _N^{-1}  (\Lambda (F_{N}) -C    \big)$ converges  in law to normal random variable. 

In the case of random variables in Wiener chaos,  from the main result in \cite{PeTu} (see Theorem \ref{m4m} later in the paper), we can prove that, under the convergence in law of $F _{N}^{(i)}$, $i=1,..,d$ to the Gaussian law, the quantities $\Lambda _{i,j} - C_{i,j}, \Lambda _{j,i} - C_{i,j}$ converges to zero in $L ^{2}(\Omega)$. Thus condition {\bf [C1\!\! ]} intuitively  means that these two terms have the same rate of convergence to zero.

\subsection{The truncation}\label{trunc}

Let us introduce the truncation sequence $(\Psi_{N}) _{N\geq 1}$ which will be widely used throughout our work.  Its use allows to avoid the condition on the existence of negative moments for the Malliavin derivative of $F_{N}$ (\ref{neg}).

We consider the function $\Psi _{N}\in \mathbb{D} ^{\ell,p} $ for  $p>1$  and $\ell \geq d+3$ such that

\begin{equation}\label{psin}
\Psi _{N}= \psi \left(  \left(  K\frac{\left| \Lambda (F_{N}) -C\right| _{d\times d} }{\gamma _{N} ^{\delta} }\right)^ {2} \right)
\end{equation}
with some $K, 0<\delta <1$, where $\psi \in C^{\infty} (\mathbb{R}; [0,1]) $ is such that $\psi (x)= 1$ is $\vert x \vert \leq \frac{1}{2}$ and $\psi(x)=0$ if $\vert  x \vert \geq 1$.  The square in the argument pf $\psi$ in order to have $\Psi _{N}$ differentiable in the Mlliavin sense.

From this definition, clearly we have
\begin{enumerate}
\item $0\leq \Psi_{N} \leq 1$  for every $N\geq 1$.

\item $\Psi _{N} \to _{N\to \infty} 1$ in probability.

\end{enumerate} 
 A crucial fact in our computations is that on the set $\{ \Psi_{N}>0\} $ we have

\begin{equation*}
\left| \Lambda (F_{N}) -C\right| _{d\times d}  \leq \frac{1}{K}\gamma _{N} ^{\delta}.
\end{equation*}
This also implies (see e.g. \cite{Yos1}) that there exists two positive constants $c_{1}, c_{2}$ such that for $N$ large enough 
\begin{equation}
\label{nd}
 c_{1} \leq\det \Lambda (F_{N})  \leq c_{2}.
\end{equation}
 
\subsection{Notation}
For every  ${\bf x}=(x_{1},..., x_{d}) \in \mathbb{R} ^{d}$ and $\boldsymbol{\alpha}=  (\alpha _{1},..., \alpha _{d}) \in \mathbb{Z} _{+} $ we use the notation 
\begin{equation*}
{\bf x} ^{\boldsymbol{\alpha} } = x_{1} ^{\alpha _{1}}....x_{d} ^{\alpha _{d}} \mbox{ and } \frac{\partial ^{\boldsymbol{\alpha}}}{\partial x^{\boldsymbol{\alpha}}}=\frac{\partial ^{\alpha_{1}}}{\partial  x_{1} ^{\alpha_{1}}}....\frac{\partial ^{\alpha_{d}}}{\partial  x_{d} ^{\alpha_{d}}}
\end{equation*}

If $ x, y \in \mathbb{R} ^{d}$, we denote by $\langle x, y\rangle _{d}:=\langle x, y\rangle $ their scalar product in $\mathbb{R} ^{d}$. By $\vert x \vert _{d}:=\vert x\vert $  we will denote the Euclidean norm of $x$.

The following will be fixed in the sequel:
\begin{itemize}
\item The sequence $(F_{N}) _{N\geq 1}$ defined in Section \ref{21}.

\item The truncation sequence $(\Psi _{N}) _{N\geq 1} $ is as in Section \ref{trunc}.

\item  $o_{\lambda} (\gamma _{N}) $ denotes a deteministic sequence that depends on $\lambda$ such that $\frac{o _{\lambda } (\gamma _{N})}{\gamma _{N}} \to_{N \to \infty }0$. By $ o(\gamma_{N})$ we refer as usual to a sequence such that $\frac{o (\gamma _{N})}{\gamma _{N}} \to_{N}0.$

\item By $\Phi (\lambda)= e ^{-\frac{\lambda ^{2}}{2}}$ we denote the characteristic function of $Z\sim N(0,1)$ and by $p(x)=\frac{1}{\sqrt{2\pi}}e ^{-\frac{x ^{2}}{2}}$ its density.

\item By $c, C, C_{1}, ...$ we denote  generic strictly positive constants that may vary from line to line and by $\langle \cdot , \cdot \rangle $ we denote the Euclidean scalar product $\mathbb{R} ^{d}$. On the other hand, we will keep the subscript for the inner product in $H$, denoted by $\langle \cdot,\cdot \rangle _{H}$. 

\item We will use bold letters to indicate vectors in $\mathbb{R}^{d}$ when we need to differentiate them from their one-dimensional components.
\end{itemize}

\section{Asymptotic behavior of  the characteristic function}

The first step in finding the asymptotic expansion for the sequence $(F_{N})_{N\geq 1}$ which satisfies  {\bf  [C1\!\! ] } - {\bf  [C3\!\! ] } is to analyze the asymptotic behavior  as $N\to \infty$ of the characteristic function of $F_{N}$. In order to avoid troubles related to the integrability over the whole real line, we will work under truncation, in the sense that we will always keep the multiplicative factor $\Psi_{N}$ given by  (\ref{psin}).  

Let $ \boldsymbol{\lambda}= (\lambda _{1},..., \lambda _{d}) \in \mathbb{R } ^{d}$ and $\theta \in [0,1]$ and let us consider the truncated interpolation

\begin{equation}\label{interpol}
\varphi _{N} ^{\Psi} (\theta, \boldsymbol{\lambda} ) = \mathbf{E} \left( \Psi _{N} e ^{{\tt i}\theta \langle \boldsymbol{\lambda} , F _{N} \rangle - (1-\theta ^{2} ) \frac{\boldsymbol{\lambda} ^{T}C\boldsymbol{\lambda} }{2}}\right).
\end{equation}
Notice that $\varphi_{N} (1, \boldsymbol{\lambda}) = \mathbf{E} (\Psi_{N} e ^{i\boldsymbol{\lambda}, F_{N} \rangle })$, the "truncated" characteristic function of $F_{N}$, while $\varphi_{N} (0, \boldsymbol{\lambda}) =(\mathbf{E} \Psi_{N})e  ^{-\frac{\boldsymbol{\lambda} ^{T}C\boldsymbol{\lambda} }{2}}$, the "truncated"characteristic function of the limit in law of $F_{N}$.  Therefore, it is useful to analyze the behavior of the derivative of $\varphi _{N} ^{\Psi} (\theta, \boldsymbol{\lambda} ) $ with respect to the variable $\theta$.
\vskip0.2cm

We have the following result.

\begin{prop}\label{p101}
Assume {\bf  [C1\!\! ] } - {\bf  [C3\!\! ] }  and suppose that

\begin{equation}
\label{z12} \mathbf{E} \left(Z_{2} ^{(k,l)} | Z _{1} \right)= \sum_{a=1} ^{d} \rho_{k,l} ^{a} Z _{1} ^{(a) }
\end{equation}
for every $k,l=1,..,d$, where $(Z_{1}, Z_{2})$ is the Gaussian random vector from {\bf  [C1\!\! ]. } Then we have the convergence

 \begin{equation}\label{conv1}
 \gamma _{N } ^{-1} \frac{\partial}{\partial \theta} \varphi _{N} ^{\Psi} (\theta, \boldsymbol{\lambda}) \to _{N\to \infty} 
-{\tt i}\theta ^{2} \left( \sum_{i, j, k, a=1}^{d} \lambda _{j} \lambda _{k}\lambda _{i} \rho_{j,k} ^{a} C_{ia}   \right) e ^{-\frac{\boldsymbol{\lambda} ^{T} C \boldsymbol{\lambda} }{2}}.
\end{equation}
\end{prop}
{\bf Proof: } By differentiating (\ref{interpol}) with respect to $\theta$  and using the formula $ F_{N}^{(i)}=\delta D (-L) ^{-1} F _{N} ^{(i)}$ for every $i=1,.., d$,  we can write 

\begin{eqnarray*}
\frac{\partial}{\partial \theta} \varphi _{N} ^{\Psi} (\theta, \boldsymbol{\lambda}) &=&\mathbf{E}\left( \Psi _{N} e ^{{\tt i}\theta \langle \boldsymbol{\lambda}, F _{N} \rangle - (1-\theta ^{2} ) \frac{\boldsymbol{\lambda} ^{T}C\boldsymbol{\lambda} }{2}} \left( i \langle \boldsymbol{\lambda} , F_{N} \rangle + \theta \boldsymbol{\lambda} ^{T} C \boldsymbol{\lambda}  \right) \right)\\
&=&{\tt i} \mathbf{E} \left( \Psi _{N} e ^{{\tt i}\theta \langle \boldsymbol{\lambda} , F _{N} \rangle - (1-\theta ^{2} ) \frac{\boldsymbol{\lambda} ^{T}C\boldsymbol{\lambda} }{2}}  \sum_{j=1} ^{d} \lambda _{j} \delta D (-L) ^{-1} F_{N} ^{(j)}\right)\\
&&+ \theta \sum_{j,k=1} ^{d} \lambda _{j} \lambda _{k} C_{j,k} \mathbf{E}\left( \Psi _{N} e ^{{\tt i}\theta \langle \boldsymbol{\lambda} , F _{N} \rangle - (1-\theta ^{2} ) \frac{\boldsymbol{\lambda} ^{T}C\boldsymbol{\lambda} }{2}} \right)
\end{eqnarray*}
and by the duality relationship (\ref{dua})
\begin{eqnarray*}
\frac{\partial}{\partial \theta} \varphi _{N} ^{\Psi} (\theta,\boldsymbol{\lambda}) &=&{\tt i} \sum_{j=1} ^{d}\lambda _{j} \mathbf{E} \left( \Psi _{N} e ^{{\tt i}\theta \langle \boldsymbol{\lambda} , F _{N} \rangle - (1-\theta ^{2} ) \frac{\boldsymbol{\lambda} ^{T}C\boldsymbol{\lambda }}{2}} {\tt i}\theta \langle D(-L) ^{-1} F_{N} ^{(j)} , D\langle \boldsymbol{\lambda}, F_{N} \rangle  \rangle _{H} \right)\\
&&+{\tt i} \sum_{j=1} ^{d} \lambda _{j} \mathbf{E} \left( e ^{{\tt i}\theta \langle \boldsymbol{\lambda} , F _{N} \rangle - (1-\theta ^{2} ) \frac{\boldsymbol{\lambda} ^{T}C\boldsymbol{\lambda }}{2}}\langle D\Psi_{N} , D(-L) ^{-1} F _{N} ^{(j)}\rangle _{H}\right)\\
&&+\theta  \sum_{j,k=1} ^{d} \lambda _{j} \lambda _{k} C_{j,k}\mathbf{E} \left( \Psi _{N} e ^{{\tt i}\theta \langle \boldsymbol{\lambda} , F _{N} \rangle - (1-\theta ^{2} ) \frac{\boldsymbol{\lambda} ^{T}C\boldsymbol{\lambda }}{2}} \right).
\end{eqnarray*}
We obtain, since $D\langle \boldsymbol{\lambda}, F_{N} \rangle= \sum_{k=1} ^{d} \lambda _{k}DF_{N} ^{(k)}$,

\begin{eqnarray}
 \frac{\partial}{\partial \theta} \varphi _{N} ^{\Psi} (\theta, \boldsymbol{\lambda} ) &=&-\theta \sum_{j,k=1} ^{d} \lambda _{j} \lambda _{k} \mathbf{E} \left( \Psi _{N} e ^{{\tt i}\theta \langle \boldsymbol{\lambda} , F _{N} \rangle - (1-\theta ^{2} ) \frac{\boldsymbol{\lambda} ^{T}C\boldsymbol{\lambda }}{2}} \left(  \langle DF_{N} ^{(k)}, D(-L) ^{-1} F _{N} ^{(j)}\rangle_{H} - C_{j,k}\right) \right) \nonumber \\
&&+{\tt i} \sum_{j=1} ^{d}\lambda _{j}  \mathbf{E} \left( e ^{{\tt i}\theta \langle \boldsymbol{\lambda} , F _{N} \rangle - (1-\theta ^{2} ) \frac{\boldsymbol{\lambda} ^{T}C\boldsymbol{\lambda }}{2}}\langle D\Psi_{N} , D(-L) ^{-1} F _{N} ^{(j)}\rangle _{H} \right)\nonumber\\
&:=& A_{N}(\theta, \boldsymbol{\lambda})+ B _{N} (\theta, \boldsymbol{\lambda}).  
\label{31o-1}
\end{eqnarray}

\noindent Using the definition of the truncation function $\Psi_{N}$, we can prove that 
\begin{equation}
\label{11m-1}
\gamma_{N} ^{-1} B_{N} (\theta, \boldsymbol{\lambda})= \gamma _{N } ^{-1}\sum_{j=1} ^{d} \mathbf{E} \left( e ^{{\tt i}\theta \langle \boldsymbol{\lambda} , F _{N} \rangle - (1-\theta ^{2} ) \frac{\boldsymbol{\lambda} ^{T}C\boldsymbol{\lambda }}{2}}\langle D\Psi_{N} , D(-L) ^{-1} F _{N} ^{(j)}\rangle _{H} \right) \to_{N\to \infty } 0.
\end{equation}
Indeed, by the chain rule (\ref{chain})
$$D\Psi_{N}= \psi ' \left(  \left(  K\frac{\left| \Lambda (F_{N}) -C\right| _{d\times d} }{\gamma _{N} ^{\delta} }\right) ^ {2}\right) D\left( \frac{\left| \Lambda (F_{N}) -C\right| _{d\times d} }{\gamma _{N} ^{\delta} }\right)^ {2}$$
and by using {\bf  [C2\!\! ]}-{\bf  [C3\!\! ]} and (\ref{psin}), we get for every $p>1$
\begin{eqnarray*}
\gamma_{N} ^{-1}\vert  B_{N} (\theta, \boldsymbol{\lambda} ) \vert &&\leq c\gamma_{N} ^{-1} \mathbf{E} \left|\langle D\Psi_{N} , D(-L) ^{-1} F _{N} ^{(j)}\rangle _{H} \right| \\
&&\leq c\gamma_{N} ^{-1}  P \left( \left| \Lambda (F_{N}) -C\right| _{d\times d} \geq \frac{1}{2} \gamma_{N} ^{\delta} \right) ^{\frac{1}{2}}\\
&&\leq c \mathbf{E} \left( \left| \Lambda (F_{N}) -C\right| _{d\times d} \right) ^{p} \gamma_{N}^{-1-\delta p}\leq c \gamma_{N}^{ p(1-\delta)- 1}
\end{eqnarray*}
 where we used again {\bf  [C2\!\! ]}. Since $\delta <1$ and $p$ is arbitrarly large, we obtain the convergence (\ref{11m-1}).

On the other hand, by assumption {\bf  [C1\!\! ]}, we have the convergence 

\begin{equation}
\label{1n-1}
\gamma_{N} ^{-1} A_{N} (\theta, \boldsymbol{\lambda} )\to_{N\to \infty} -\theta\sum_{j,k=1} ^{d} \lambda _{j} \lambda _{k}\mathbf{E}\left( e ^{{\tt i}\theta \langle \boldsymbol{\lambda} , Z_{1} \rangle - (1-\theta ^{2})\frac{\boldsymbol{\lambda}^{T} C \boldsymbol{\lambda}}{2}}Z _{2} ^{(k,l)}\right).
\end{equation}
Using our assumption (\ref{z12}), we can express the above limit in a more explicit way. We have 
\begin{eqnarray}
\mathbf{E} \left( e ^{{\tt i}\theta \langle \boldsymbol{\lambda} , Z_{1} \rangle - (1-\theta ^{2})\frac{\boldsymbol{\lambda} ^{T} C \boldsymbol{\lambda} }{2}}Z _{2} ^{(j,k)}\right)
&=&\mathbf{E} \left( e ^{{\tt i}\theta \langle  \boldsymbol{\lambda}  , Z_{1} \rangle - (1-\theta ^{2})\frac{\boldsymbol{\lambda} ^{T} C \boldsymbol{\lambda} }{2}} \mathbf{E}(Z_{2} ^{(j, k) }| Z_{1})\right)\nonumber\\
&=&\mathbf{E} \left( e ^{{\tt i}\theta \langle \boldsymbol{\lambda} , Z_{1} \rangle - (1-\theta ^{2})\frac{\boldsymbol{\lambda} ^{T} C\boldsymbol{\lambda} }{2}}\sum_{a=1} ^{d} \rho _{j, k} ^{a} Z _{1} ^{(a) }\right)\label{1n-2}
\end{eqnarray}
and, moreover, for every $a=1,..,d$, from the differential equation verified by the characteristic function of the normal distribution, 

\begin{eqnarray}
&&\mathbf{E} \left( e ^{{\tt i}\theta \langle \boldsymbol{\lambda} , Z_{1} \rangle - (1-\theta ^{2})\frac{\boldsymbol{\lambda} ^{T} C \boldsymbol{\lambda} }{2}} Z _{1} ^{(a) }\right) =e ^{- (1-\theta ^{2})\frac{\boldsymbol{\lambda} ^{T} C \boldsymbol{\lambda} }{2}} \mathbf{E}\left( e ^{{\tt i}\theta \langle\boldsymbol{\lambda} , Z_{1} \rangle}Z _{1} ^{(a) }\right) \nonumber\\
&=&e ^{- (1-\theta ^{2})\frac{\boldsymbol{\lambda}^{T} C\boldsymbol{\lambda} }{2}} \frac{-1}{{\tt i}\theta} (\sum_{i=1} ^{d} \lambda _{i}C_{ia})\theta ^{2}  e ^{-\theta ^{2}\frac{\boldsymbol{\lambda} ^{T} C \boldsymbol{\lambda} }{2}}\nonumber \\
&=&e ^{-\frac{\boldsymbol{\lambda} ^{T} C \boldsymbol{\lambda} }{2}} \frac{-\theta}{{\tt i}}  (\sum_{i=1} ^{d} \lambda _{i}C_{ia}) ={\tt i}\theta e ^{-\frac{\boldsymbol{\lambda}^{T} C \boldsymbol{\lambda} }{2}}  \sum_{i=1} ^{d} \lambda _{i}C_{ia}.\label{1n-3}
\end{eqnarray}
Thus, by replacing (\ref{1n-2}) and (\ref{1n-3}) in (\ref{1n-1}), we obtain
\begin{equation}
\label{31o-2}
\gamma_{N} ^{-1} A_{N} (\theta, \boldsymbol{\lambda} )\to_{N\to \infty}-{\tt i}\theta ^{2} \left( \sum_{i, j, k, a=1}^{d} \lambda _{j} \lambda _{k}\lambda _{i} \rho_{j,k} ^{a} C_{ia}   \right) e ^{-\frac{\boldsymbol{\lambda}^{T} C \boldsymbol{\lambda} }{2}}.
\end{equation}
Relation (\ref{31o-2}), together with with (\ref{11m-1}) and (\ref{31o-1}) will imply the conclusion. \qed

\vskip0.2cm

We will also need to analyze the asymptotic behavior  of the partial derivatives with respect to the variable $\boldsymbol{\lambda}$ of (\ref{interpol}). Let us first introduce some notation borrowed from \cite{Yos3}.

For every $u, z \in \mathbb{R} ^{d},  \boldsymbol{\alpha} \in \mathbb{Z}_{+} ^{d}$ and for every symmetric $d\times d$ matrix $C$, we will  denote by 
\begin{equation}\label{pa}
P_{\boldsymbol{\alpha}}(u,z, C)= e ^{-{\tt i}\langle u, z \rangle _{d} - \frac{1}{2} u^{T} Cu} (-{\tt i}) ^{\vert \boldsymbol{\alpha}\vert }  \frac{\partial ^{\boldsymbol{\alpha}}}{\partial u ^{\boldsymbol{\alpha}}} e ^{{\tt i}\langle u, z \rangle _{d} +\frac{1}{2} u^{T} Cu} 
\end{equation}
 where  $\vert \boldsymbol{\alpha}\vert := \alpha_{1}+...+\alpha_{d}$  if $\boldsymbol{\alpha}=(\alpha_{1},.., \alpha_{d})\in \mathbb{Z}_{+} ^{d}$. Then clearly, 
\begin{equation}
\label{31o-4}
(-{\tt i}) ^{\vert \boldsymbol{\alpha}\vert }  \frac{\partial ^{\boldsymbol{\alpha}}}{\partial u ^{\boldsymbol{\alpha}}} e ^{{\tt i}\langle u, z \rangle _{d} +\frac{1}{2} u^{T} Cu}=e ^{{\tt i}\langle u, z \rangle _{d} +\frac{1}{2} u^{T} Cu}P_{\boldsymbol{\alpha}} (u,z, C).
\end{equation}

In the next lemma we recall the explicit expression of the polynomial $P_{\boldsymbol{\alpha}}$ together on some estimates on this polynomial. We refer to Lemma 1 in \cite{Yos3} for the proof.

\begin{lemma}\label{l-101} If $P_{\boldsymbol{\alpha}}$ is given by (\ref{pa}), then for every $u, z \in \mathbb{R} ^{d}$ and evey $d\times d$ matrix $C$, we have
$$ P_{\boldsymbol{\alpha}} (u,z, C)=\sum_{0\leq \beta \leq \boldsymbol{\alpha}} C_{\boldsymbol{\alpha}} ^{\beta} (-{\tt i}) ^{\beta} z ^{\boldsymbol{\alpha} -\beta}   S_{\beta} (u,C)$$
with
$$S_{\beta} (u, C)= e ^{-\frac{1}{2} u^{T} Cu} \frac{\partial ^{\beta}}{\partial u ^{\beta} } e ^{\frac{1}{2}
u^{T}Cu}.$$
Moreover, we have the estimate
$$\vert S_{\beta} (u, C) \vert \leq   \sum_{j=0} ^{\vert \beta \vert} c_{j} ^{\vert\beta \vert } \vert u \vert ^{j} \vert C \vert ^{( j+\vert \beta \vert )/2}$$
for  every $u \in \mathbb{R} ^{d}$ and every $d\times d$ matrix $C$.
\end{lemma}

Let us now state and prove the asymptotic behavior of the partial derivatives of the truncated interpolation $\varphi ^{\Psi}_{N}$.

\begin{prop}
Assume {\bf  [C1\!\! ] }-{\bf  [C3\!\! ] }  and (\ref{z12}). Let $\varphi ^{\Psi} _{N}$ be given by (\ref{interpol}). Then for every $\boldsymbol{\alpha}\in \mathbb{Z}_{+} ^ {d}$, we have
\begin{equation*}
\gamma_{N} ^{-1} \frac{ \partial ^{\boldsymbol{\alpha}}}{ \partial \boldsymbol{\lambda}^{\boldsymbol{\alpha}}} \frac{\partial}{\partial \theta} \varphi _{N} ^{\Psi} (\theta, \boldsymbol{\lambda}) \to _{N\to \infty} \frac{ \partial ^{\boldsymbol{\alpha}}}{ \partial \boldsymbol{\lambda}^{\boldsymbol{\alpha}}}\left[ 
-{\tt i}\theta ^{2} \left( \sum_{i, j, k, a=1}^{d} \lambda _{j} \lambda _{k}\lambda _{i} \rho_{j,k} ^{a} C_{ia}   \right) e ^{-\frac{\boldsymbol{\lambda }^{T} C \boldsymbol{\lambda} }{2}}\right].
\end{equation*}
\end{prop}
{\bf Proof: }If $\boldsymbol{\alpha}= (\alpha_{1},  ..., \alpha_{d}) \in \mathbb{Z}_{+}^{d} $, then by (\ref{31o-1}) we have
 \begin{equation*}
\frac{ \partial ^{\boldsymbol{\alpha}}}{ \partial \boldsymbol{\lambda}^{\boldsymbol{\alpha}}} \frac{\partial}{\partial \theta} \varphi _{N} ^{\Psi} (\theta, \boldsymbol{\lambda}) = \frac{ \partial ^{\boldsymbol{\alpha}}}{ \partial \boldsymbol{\lambda}^{\boldsymbol{\alpha}}} A_{N}(\theta, \boldsymbol{\lambda} )+ \frac{ \partial ^{\boldsymbol{\alpha}}}{ \partial \boldsymbol{\lambda}^{\boldsymbol{\alpha}}}B_{N} (\theta, \boldsymbol{\lambda})
\end{equation*}
with $A_{N}, B_{N}$ defined in (\ref{31o-1}). Now, since by (\ref{31o-4})
$$\frac{\partial ^{\boldsymbol{\alpha}} }{\partial \boldsymbol{ \lambda} ^{\alpha}}  e ^{{\tt i}\theta \langle \boldsymbol{\lambda} , F _{N} \rangle - (1-\theta ^{2} ) \frac{\boldsymbol{\lambda} ^{T}C\boldsymbol{\lambda }}{2}}= {\tt i} ^{\vert \boldsymbol{\alpha}\vert }P_{\boldsymbol{\alpha}} (\boldsymbol{\lambda}, \theta F_{N} , -(1-\theta ^{2}) C)e ^{{\tt i}\theta \langle \boldsymbol{\lambda} , F _{N} \rangle - (1-\theta ^{2} ) \frac{\boldsymbol{\lambda} ^{T}C\boldsymbol{\lambda }}{2}}$$ 
we can write

\begin{eqnarray}
&& \frac{\partial ^{\boldsymbol{\alpha}} }{\partial\boldsymbol{\lambda} ^{\alpha}} \frac{\partial}{\partial \theta} \varphi _{N} ^{\Psi} (\theta, \boldsymbol{\lambda})\nonumber \\
&=&-\theta \sum_{j,k=1} ^{d} \lambda_{j}\lambda _{k}\mathbf{E} \left( \Psi _{N} e ^{{\tt i}\theta \langle \boldsymbol{\lambda} , F _{N} \rangle - (1-\theta ^{2} ) \frac{\boldsymbol{\lambda} ^{T}C\boldsymbol{\lambda }}{2}} {\tt i} ^{\boldsymbol{\alpha}} P_{\boldsymbol{\alpha}}  (\boldsymbol{\lambda}, \theta F_{N}, -(1-\theta ^{2}),C) \left(  \langle DF_{N} ^{(k)}, D(-L) ^{-1} F _{N} ^{(j)}\rangle _{H}- C_{j,k}\right) \right)\nonumber  \\
&&+{\tt i} \sum_{j=1} ^{d} \lambda _{j}\mathbf{E} \left( \Psi _{N} e ^{{\tt i}\theta \langle \boldsymbol{\lambda} , F _{N} \rangle - (1-\theta ^{2} ) \frac{\boldsymbol{\lambda} ^{T}C\boldsymbol{\lambda }}{2}}  {\tt i} ^{\boldsymbol{\alpha}} P_{\boldsymbol{\alpha}}  ( \boldsymbol{\lambda}, \theta F_{N}, -(1-\theta ^{2}),C)  \langle D\Psi_{N} , D(-L) ^{-1} F _{N} ^{(j)}\rangle _{H}\right)\label{31o-5}
\end{eqnarray}

As in the proof of (\ref{11m-1}), the second summand above multiplied by $\gamma_{N} ^{-1}$  converges to zero as $N\to \infty$.  Concerning the first summand in (\ref{31o-5}), the hypothesis {\bf  [C1\!\! ]} implies that it converges as $N\to \infty$, to 
 
\begin{eqnarray*}
&&-\theta\sum_{j,k=1} ^{d} \lambda _{j} \lambda _{k}\mathbf{E} \left( {\tt i} ^{\boldsymbol{\alpha}} P_{\boldsymbol{\alpha}}  (\boldsymbol{\lambda}, \theta Z_{1}, -(1-\theta ^{2})C) e ^{{\tt i}\theta \langle\boldsymbol{\lambda} , Z_{1} \rangle - (1-\theta ^{2})\frac{\boldsymbol{\lambda}^{T} C \boldsymbol{\lambda} }{2}}Z _{2} ^{(k,l)}\right) \\
&=& -\theta\sum_{j,k=1} ^{d} \lambda _{j} \lambda _{k}\mathbf{E}\left( \left( \frac{\partial ^{\boldsymbol{\alpha}} }{\partial \boldsymbol{\lambda} ^{\alpha}} e ^{{\tt i}\theta \langle \boldsymbol{\lambda}, Z_{1} \rangle - (1-\theta ^{2})\frac{\boldsymbol{\lambda} ^{T} C \boldsymbol{\lambda} }{2}} \right)Z _{2} ^{(k,l)}\right) \\
&=&\frac{ \partial ^{\boldsymbol{\alpha}}}{ \partial \boldsymbol{\lambda}^{\boldsymbol{\alpha}}}\left[ 
-{\tt i}\theta ^{2} \left( \sum_{i, j, k, a=1}^{d}\lambda _{j} \lambda _{k}\lambda _{i} \rho_{j,k} ^{a} C_{ia}   \right) e ^{-\frac{\boldsymbol{\lambda }^{T} C \boldsymbol{\lambda} }{2}}\right].
\end{eqnarray*} 
Then the conclusion is obtained. \qed

\begin{remark}
If $d=1$, relation (\ref{conv1}) becomes
\begin{equation*}
\gamma _{N } ^{-1} \frac{\partial}{\partial \theta} \varphi _{N} ^{\Psi} (\theta, \lambda) \to _{N\to \infty} 
-{\tt i} \theta ^{2}\lambda ^{3} \rho e ^{-\frac{\lambda ^{2}}{2}}
\end{equation*}
where $\rho=\mathbf{E}(Z_{2}/Z_{1})$.
\end{remark}

The following integration by parts formula plays an important role in the sequel. Its proof is a slightly adaptation of the proof of Proposition 1 in \cite{Yos3}. For its proof the nondegeneracy (\ref{nd}) is needed.

\begin{lemma}\label{key} 
Assume {\bf  [C2\!\! ]} -{\bf  [C3\!\! ]}.  Consider a function $f \in C _{p}^{l} (\mathbb{R} ^{d})$ with $l\geq 0$ and let $G$ be a random variable in $ \mathbb{D} ^{l, \infty}$. Then it holds
\begin{equation}\label{ibp}
\mathbf{E} \left( \frac{\partial ^{ l}}{\partial x ^{l}} f(F_{N}) G \Psi _{N}\right) =\mathbf{E}\left( f(F_{N}) H_{l} (G\Psi_{N} ) \right)
\end{equation}
where the random variable $H_{l} (G\Psi_{N} ) $ belongs to $ L^{p} (\Omega)$ for every $p\geq 1$.
\end{lemma}

We will use the integration by parts formula in the above Lemma \ref{key} to show that $\frac{\partial ^{\boldsymbol{\alpha}} }{\partial\boldsymbol{\lambda}  ^{\boldsymbol{\alpha}}} \frac{\partial}{\partial \theta} \varphi _{N} ^{\Psi} (\theta, \boldsymbol{\lambda})$ is dominated by an integrable function uniformly with respect to $(\theta, \boldsymbol{\lambda})\in [0,1] \times \mathbb{R} ^{d}.$

\begin{lemma}\label{170223-4} Assume {\bf  [C2\!\! ]}-{\bf  [C3\!\! ]}. 
For every $\boldsymbol{\alpha}\in\bbZ_{+} ^{d}$, there exists a positive constant $C_{\boldsymbol{\alpha}}$ such that 
\begin{equation}\label{310-5} 
\gamma _N^{-1}\bigg|\partial_ {\boldsymbol{\lambda}}^{\boldsymbol{\alpha}} \frac{\partial }{\partial \theta}\varphi_{F_{N}}  ^{\Psi} (\theta,\boldsymbol{\lambda})\bigg|
\leq
C_{\boldsymbol{\alpha}} (1+|\boldsymbol{\lambda}|)^{-\ell+2}
\end{equation}
for all $\boldsymbol{\lambda}\in\bbR ^{d} $, $\theta\in[0,1]$ and $N\geq 1$. 
\end{lemma}
{\bf Proof: }
First we consider the case $\boldsymbol{\alpha}=0\in \mathbb{Z} _{+}^{d}$. Consider the function $f(x)= e ^{{\tt i}\theta \langle \boldsymbol{\lambda} , x\rangle - (1-\theta ^{2}) \frac{\boldsymbol{\lambda} ^{T} C\boldsymbol{\lambda} }{2}}$ for $x\in \mathbb{R} ^{d}$ with  
$\frac{ \partial ^{\boldsymbol{\alpha }}}{\partial x ^{\boldsymbol{\alpha}}}f(x)= f(x) (i\boldsymbol{\lambda} ) ^{\boldsymbol{\alpha}}.$

For $\theta\in(0,1]$, by using  (\ref{31o-1}) and the integration by parts formula (\ref{ibp}) to the function $f$ $\ell$-times, we obtain

\begin{eqnarray*}
&&({\tt i}\theta\boldsymbol{\lambda})^\ell
\frac{\partial }{\partial \theta}\varphi_{F_{N}}  ^{\Psi} (\theta,\boldsymbol{\lambda})\\
&=&
 -\theta ({\tt i}\theta\boldsymbol{\lambda})^\ell\sum_{j,k=1} ^{d}\lambda _{j}\lambda _{k} \mathbf{E} \left( \Psi _{N} e ^{{\tt i}\theta \langle \boldsymbol{\lambda} , F _{N} \rangle - (1-\theta ^{2} ) \frac{\boldsymbol{\lambda} ^{T}C\boldsymbol{\lambda }}{2}} \left(  \langle DF_{N} ^{(k)}, D(-L) ^{-1} F _{N} ^{(j)}\rangle _{H}- C_{j,k}\right) \right) \\
&&+{\tt i} ({\tt i}\theta\boldsymbol{\lambda})^\ell\sum_{j=1} ^{d} \lambda _{j} \mathbf{E} \left( e ^{{\tt i}\theta \langle \boldsymbol{\lambda} , F _{N} \rangle - (1-\theta ^{2} ) \frac{\boldsymbol{\lambda} ^{T}C\boldsymbol{\lambda }}{2}}\langle D\Psi_{N} , D(-L) ^{-1} F _{N} ^{(j)}\rangle _{H} \right)\\
&=& -\theta \sum_{j,k=1} ^{d} \lambda _{j}\lambda _{k}\mathbf{E} \left( \Psi _{N} e ^{{\tt i}\theta \langle \boldsymbol{\lambda} , F _{N} \rangle - (1-\theta ^{2} ) \frac{\boldsymbol{\lambda} ^{T}C\boldsymbol{\lambda }}{2}} H_\ell\bigg(\big(q^{-1}\langle DF_{N},DF_{N}\rangle _{H}-1\big)\Psi_{N}\bigg) \right)\\
&&+\sum_{j=1} ^{d}\lambda _{j} \mathbf{E} \left( e ^{{\tt i}\theta \langle \boldsymbol{\lambda} , F _{N} \rangle - (1-\theta ^{2} ) \frac{\boldsymbol{\lambda} ^{T}C\boldsymbol{\lambda }}{2}}H_\ell\bigg(\langle i\boldsymbol{\lambda} (-q)^{-1}DF_{N},D\Psi_{N}\rangle _{H}\bigg)\right)
\end{eqnarray*}

We use this estimate for $\theta\geq1/2$. 
For $\theta\in[0,1/2]$, the inequality is trivial thanks to the exponential function 
$\exp(-(1-\theta^2) \frac{\boldsymbol{\lambda} ^{T}C\boldsymbol{\lambda}}{2})$. 
Thus we obtained the desired estimate when $\boldsymbol{\alpha}=0$. 
For $\boldsymbol{\alpha}\in \mathbb{Z}_{+} ^{d}, \boldsymbol{\alpha} \not=(0,...0)$, we can follow a similar argument with the estimate 
\beas 
\sup_{\theta,\boldsymbol{\lambda}}
\big((1-\theta^2)(1+\boldsymbol{\lambda}^2)\big)^m\exp\big(-(1-\theta^2) \frac{\boldsymbol{\lambda} ^{T}C\boldsymbol{\lambda} }{2}\big)
&<&
\infty
\eeas
for every $m\in\bbZ_+$. 
Thus we obtained the result. \qed

Let us denote by $\varphi_{N}$ the sum of the characteristic function of the law $N(0, C)$ and of the integral over $[0,1]$ with respect to $\theta$ of the limit in (\ref{conv1}), i.e.
\begin{equation}\label{fin}
\varphi_{N} (\boldsymbol{\lambda})= e ^{-\frac{ \boldsymbol{\lambda} ^{T}C\boldsymbol{\lambda} }{2}}-{\tt i}\gamma_{N} \frac{1}{3} \left( \sum_{i, j, k, a=1}^{d}\lambda _{j} \lambda _{k}\lambda _{i} \rho_{j,k} ^{a} C_{ia}   \right) e ^{-\frac{\boldsymbol{\lambda} ^{T} C\boldsymbol{\lambda} }{2}}.
\end{equation}

 The next result shows that the difference between the  truncated characteristic function of $F_{N}$ and  the characteristic function of the weak limit of the sequence $F_{N}$   can be approximated by $\varphi_{N}$. 

\begin{prop}\label{l-104} Suppose that  {\bf  [C1\!\! ]}-{\bf  [C3\!\! ]} are fulfilled for  $\ell=d+3$.  For every $N\geq 1$, let $\varphi_{N}$ be given by (\ref{fin}). Then 
\beas
\gamma_{N}^{-1} 
\int  _{\mathbb{R} ^{d}}\bigg|
\partial_{\boldsymbol{\lambda}}^{\boldsymbol{\alpha}} \bigg(
\mathbf{E}\big[e^{{\tt i}\langle \boldsymbol{\lambda},F_N}\rangle\Psi_N\big]
- \mathbf{E}[\Psi_N]\varphi_{N} (\boldsymbol{\lambda}) \bigg) \bigg| d\boldsymbol{\lambda}
&\to _{N \to \infty}&
0
\eeas
for every $\boldsymbol{\alpha}\in\bbZ_{+} ^{d} $. 
\end{prop}
{\bf Proof: }If $\boldsymbol{\alpha}\in\bbZ_{+} ^{d} $, using 
$$\varphi  _{N} ^{\Psi} (1, \boldsymbol{\lambda})= \mathbf{E}\big[e^{{\tt i}\langle \boldsymbol{\lambda},  F_N\rangle}\Psi_N\big] \mbox{ and } \varphi  _{N} ^{\Psi} (0, \boldsymbol{\lambda})= e ^{-\frac{\boldsymbol{\lambda} ^{T} C \boldsymbol{\lambda} }{2}}$$
 we can write
\begin{eqnarray*}
&&
\int_{\bbR ^{d}} \bigg|\gamma_{N}^{-1}
\partial_{\boldsymbol{\lambda}}^{\boldsymbol{\alpha}} \bigg(
\mathbf{E}\big[e^{{\tt i}\langle \boldsymbol{\lambda}, F_N\rangle}\Psi_N\big]
-\mathbf{E}[\Psi_N] \varphi_{N} (\boldsymbol{\lambda})\bigg)\bigg|d\boldsymbol{\lambda}
\\
&=&
\int_\bbR \bigg|\partial_{\boldsymbol{\lambda}}^{\boldsymbol{\alpha}} 
\int_0^1\bigg(
\gamma_{N}^{-1}\partial_\theta\varphi^\Psi_N(\theta,\boldsymbol{\lambda})
-(-{\tt i})\theta ^{2} \mathbf{E}[\Psi_N] \left( \sum_{i, j, k, a=1}^{d} \lambda _{j} \lambda _{k}\lambda _{i} \rho_{j,k} ^{a} C_{ia}   \right) e ^{-\frac{\boldsymbol{\lambda }^{T} C \boldsymbol{\lambda} }{2}}\bigg)d\theta\bigg|d\boldsymbol{\lambda}
\\
&\leq&\int_{ \bbR ^{d}}\int_0^1\bigg|
\gamma_{N}^{-1}\partial_{\boldsymbol{\lambda}}^{\boldsymbol{\alpha}} \partial_\theta\varphi^\Psi_N(\theta,\boldsymbol{\lambda})
-\partial_{\boldsymbol{\lambda}} ^{\boldsymbol{\alpha}}\big((-{\tt i})\theta ^{2}\left( \sum_{i, j, k, a=1}^{d}\lambda _{j} \lambda _{k}\lambda _{i} \rho_{j,k} ^{a} C_{ia}   \right) e ^{-\frac{\boldsymbol{\lambda }^{T} C \boldsymbol{\lambda} }{2}}\big)\bigg|d\theta d\boldsymbol{\lambda}\\
&&+ C \left(  \Vert 1-\Psi_{N}\Vert _{2} ^{2} \right) ^{\frac{1}{2}}.
\end{eqnarray*}
The last quantity converges to zero as $N\to\infty$ by the choice of the truncation function and also Lemma \ref{170223-4} and the dominated convergence theorem since the integrand is bounded uniformly in $(\theta,N)$. 
\qed 
\vskip0.2cm

Let $ \varphi_{F_{N}} $ be the characteristic function of the random vector $F_{N}$, i.e.  $\varphi_{F_{N}}(\boldsymbol{\lambda}) = \mathbf{E} e ^{{\tt i}\langle \boldsymbol{\lambda}, F_{N} \rangle }$ for $\boldsymbol{\lambda} \in \mathbb{R} ^{d}.$ If $\varphi_{N}$ is given by (\ref{fin}), we have 
\begin{eqnarray*}
&& \varphi_{F_{N}}(\lambda)\\
&=& \varphi_{N} (\boldsymbol{\lambda}) \mathbf{E}\Psi_{N} 
+ \mathbf{E}\left(\Psi_{N} e ^{{\tt i}\langle\boldsymbol{\lambda} , F_{N}\rangle}\right) - \mathbf{E}(\Psi_{N}) \varphi_{N}(\boldsymbol{\lambda})   +\mathbf{E} \left((1-\Psi_{N})e ^{{\tt i}\langle \boldsymbol{\lambda}, F_{N}\rangle}\right)\\
&=& \varphi_{N} (\boldsymbol{\lambda}) + \varphi_{N} (\boldsymbol{\lambda}) (\mathbf{E}\Psi_{N}  -1)+ \mathbf{E}\left(\Psi_{N} e ^{{\tt i}\langle\boldsymbol{\lambda} , F_{N}\rangle}\right) - \mathbf{E}(\Psi_{N}) \varphi_{N}(\boldsymbol{\lambda})   +\mathbf{E} \left((1-\Psi_{N})e ^{{\tt i}\langle \boldsymbol{\lambda}, F_{N}\rangle}\right)
\end{eqnarray*}

 Proposition \ref{l-104}  shows that the difference $\mathbf{E} \left(\Psi_{N} e ^{{\tt i}\langle \boldsymbol{\lambda}, F_{N}\rangle}\right) - \mathbf{E}(\Psi_{N}) \varphi_{N}(\boldsymbol{\lambda})  $ while from the definition of the truncation function $\Psi_{N}$ we see that the term $\mathbf{E} ((1-\Psi_{N})e ^{{\tt i}\langle\boldsymbol{\lambda}, F_{N}\rangle}+ \varphi_{N} (\boldsymbol{\lambda}) (\mathbf{E}\Psi_{N}  -1)$ is also small. Consequently, $\varphi_{F_{N}}(\boldsymbol{\lambda})$ can be approximated by $\varphi_{N}$.  Actually, we have 
\begin{equation}
\label{31o-6}
\varphi_{F_{N}}(\boldsymbol{\lambda}) = \varphi_{N} (\boldsymbol{\lambda}) +s_{N} (\boldsymbol{\lambda}) 
\end{equation}
where $s_{N}(\boldsymbol{\lambda})$ is a "small term'"  such that $\gamma_{N} ^{-1} s_{N} (\boldsymbol{\lambda})$ and its derivatives are dominated by the right-hand side of (\ref{31o-5}) and it that satisfies
\begin{equation}
\label{31o-7}
s_{N} (\boldsymbol{\lambda}) \leq C \Vert 1-\Psi_{N} \Vert _{p}+o_{\boldsymbol{\lambda} }(\gamma_{N})
\end{equation}
for every $p\geq1$, where   $o_{\boldsymbol{\lambda}} (\gamma _{N}) $ denotes a deteministic sequence that depends on $\lambda$ such that $\frac{o _{\boldsymbol{\lambda} } (\gamma _{N})}{\gamma _{N}} \to_{N\to \infty }0$.

 Therefore $\varphi_{N}$ is the dominant part in the asymptotic expansion of  $ \varphi_{F_{N}} $. By inverting $\varphi_{N}$, we get the approximate density.

\section{The approximate density}

In this paragraph, our purpose is to find the first and second-order term in the asymptotic expansion of the density of $F_{N}$. We will assume throughout this section  {\bf  [C1\!\! ]}-{\bf  [C3\!\! ]} are fulfilled for  $\ell=d+3$.  Let us denote by $\varphi_{F_{N}} ^{\Psi}$ the truncated characteristic function of $F_{N}$, i.e.
\begin{equation}
\label{pfn}
\varphi ^{\Psi} _{F_{N}} (\boldsymbol{\lambda})= \mathbf{E} \left(\Psi_{N} e ^{{\tt i}\langle \boldsymbol{\lambda} , F_{N}\rangle}\right)= \varphi _{N}  ^{\Psi} (1,\boldsymbol{\lambda})
\end{equation}
with $\varphi _{N} ^{\Psi}(1, \boldsymbol{\lambda} ) $ from (\ref{interpol}).

We  define the approximate density of $F_{N}$ via the Fourier inversion of the dominant part of $ p^{\Psi} _{F_{N}}$ defined by (\ref{pfn}).
\begin{definition}
For every $N\geq 1$, we define the approximate density $p_{N}$ by
\begin{equation}
\label{pn}
p_{N} (x)= \frac{1}{(2\pi) ^{d} } \int_{\mathbb{R}}^{d} e ^{-{\tt i} \langle \boldsymbol{\lambda},  x\rangle} \varphi _{N} (\boldsymbol{\lambda} ) d\boldsymbol{\lambda}, \hskip0.5cm x\in \mathbb{R}^{d}
\end{equation}
where $\varphi _{N}$ is given by (\ref{fin}).
\end{definition}

Let us first notice that in  the case $d=1$, we can obtain an explicit expression for $p_{N}$ in the next result.  Note that for $d=1$, assuming $C_{1,1}=1$ and $\mathbf{E} Z_{1} Z_{2} =\rho$,  we have from (\ref{fin})
\begin{equation}
\label{fin1}
\varphi_{N}(\lambda)= e ^ {-\frac{\lambda ^ {2}}{2}}-{\tt i} \frac{1}{3} \lambda ^{3}\rho e ^{-\frac{\lambda ^{2}}{2}}, \mbox{ for every } \lambda \in \mathbb{R}.
\end{equation}
Recall that 
$p(x)= \frac{1}{\sqrt{2\pi}} e ^{-\frac{x^{2}}{2}} $ denotes the density of the standard normal law. By $H_{k}$ we denote the Hermite polynomial of degree $k\geq 0$ given by $H_{0} (x)=1$ and for $k\geq 1$
$$H_{k}(x)= (-1) ^{k} e ^{\frac{x^{2}}{2}} \frac{d^{k} }{dx ^{k}} \left( e ^{-\frac{x^{2}}{2}} \right). $$

\begin{prop}\label{17d-2}
If $p_{N}$ is given by (\ref{pn}) then for every $x\in \mathbb{R}$, 
\begin{equation*} p_{N}(x)= p(x)-\frac{\rho}{3} \gamma _{N} H_{3}(x) p(x).
\end{equation*}
\end{prop}
{\bf Proof: } This follows from (\ref{fin1}) and the formula
\begin{equation}\label{17d-3}
\frac{1}{2\pi} \int_{\mathbb{R}}e ^{-{\tt i}\lambda x} ({\tt i} \lambda ) ^{k} \Phi (\lambda) d\lambda  = H_{k}(x) p(x)  
\end{equation}
for every $k\geq 0$ integer. Recall the notation  $ \Phi (\lambda):= e ^{-\frac{\lambda ^{2}}{2}}.$\qed

We prefer to work with the truncated density of $F_{N}$, which can be easier controlled.
\begin{definition}
The local  (or truncated) density of the random variable $F_{N}$ is defined as the inverse Fourier transform of the truncated characteristic function 

\begin{equation*}
p_{F_{N}} ^{\Psi} (x)= \frac{1}{(2\pi) ^{d}} \int_{\mathbb{R} ^{d}} e ^{-{\tt i}\langle \boldsymbol{\lambda}, x\rangle} \varphi ^{\Psi} _{F_{N}} (\boldsymbol{\lambda}) d\boldsymbol{\lambda}
\end{equation*}
for $x\in \mathbb{R} ^{d} $. 
\end{definition}

The local density is well-defined since obviously the truncated characteristic function $\varphi ^{\Psi} _{F_{N}} $ is integrable over $\mathbb{R}^{d} $.

Our purpose is to show that the local density $p_{F_{N}} ^{\Psi}$ is well-approximated by the approximate density $p_{N}$ given by (\ref{pn}) in the sense that for every $\boldsymbol{\alpha} \in \mathbb{Z}_{+} ^{d}$ 
\begin{equation*}
\sup_{x\in \mathbb{R}^{d}}\left| x^{\boldsymbol{\alpha}} \left( p_{F_{N}}^{\Psi}(x)- p_{N}(x) \right) \right| .
\end{equation*}
is of order $ o (\gamma _{N} )$. This will be showed in the next result, based on the asymptotic behavior of the characteristic function of $F_{N}$.

\begin{theorem}\label{tp}
For every $p\geq 1$,  and for every integer $\boldsymbol{\alpha} = (\alpha_{1}, .., \alpha_{d})\in \mathbb{Z }_{+} ^{d} $
\begin{equation}\label{p-p}
\sup_{x\in \mathbb{R}  ^{d}}\left| x^{\boldsymbol{\alpha}} \left( p ^{\Psi}_{F_{N}} (x) - p_{N}(x) \right) \right| \leq c  \Vert 1-\Psi_{N} \Vert _{p}  + o (\gamma_{N})= o (\gamma_{N}). 
\end{equation}
\end{theorem}
{\bf Proof: } We have, by integrating by parts 
\begin{eqnarray*}
x^{\boldsymbol{\alpha}}  \left( p_{F_{N}}^{\Psi}(x)- p_{N}(x) \right)  &=& x^{\boldsymbol{\alpha}}   \frac{1}{(2\pi) ^{d}} \int_{\mathbb{R}} e ^{-{\tt i} \langle\boldsymbol{\lambda}, x\rangle} \left( \varphi ^{\psi } _{F_{N} }(\boldsymbol{\lambda} -\varphi_{N} (\lambda) \right)d\boldsymbol{\lambda}\\
&=& \frac{1}{(2\pi ) ^{d} (-{\tt i}) ^{\vert \boldsymbol{\alpha}\vert} } \int_{\mathbb{R} ^{d}}e ^{-{\tt i} \langle\boldsymbol{\lambda}, x\rangle} \frac{ \partial ^{\boldsymbol{\alpha}}}{ \partial \boldsymbol{\lambda}^{\boldsymbol{\alpha}}} \left( \varphi ^{\psi } _{F_{N} }(\boldsymbol{\lambda}) -\varphi_{N} (\boldsymbol{\lambda}) \right)d\boldsymbol{\lambda}
\end{eqnarray*}
and then, by using (\ref{31o-6}) and Lemma \ref{170223-4} we obtain
\begin{eqnarray*}
\sup_{x\in \mathbb{R}^{d} }\left| x^{\boldsymbol{\alpha}} \left( p_{F_{N}}^{\Psi}(x)- p_{N}(x) \right) \right| &\leq &  \frac{1}{(2\pi) ^{d}}
 \int_{\mathbb{R} ^{d}} \left| \frac{d ^{\boldsymbol{\alpha}}}{d\boldsymbol{\lambda} ^{\alpha}} \left( \varphi ^{\psi } _{F_{N} }(\lambda) -\varphi_{N} (\boldsymbol{\lambda}) \right) \right| d\boldsymbol{\lambda} \\
&\leq & c \int_{\mathbb{R} ^{d}} (1+|\boldsymbol{\lambda}|)^{-\ell+2}d\boldsymbol{\lambda}\left( o (\gamma_{N}) + c\Vert 1-\Psi_{N} \Vert _{p}\right) \\
&\leq & o (\gamma_{N}) + c\Vert 1-\Psi_{N} \Vert _{p}= o (\gamma_{N})
\end{eqnarray*}
where we use $\ell-2>d$ since $\ell\geq d+3 $ in {\bf  [C2\!\! ]}. The fact that $c\Vert 1-\Psi_{N} \Vert _{p} = o (\gamma_{N})$ follows from the choice of our truncation function $\Psi_{N}$.
\qed

Let us make some comments on the above result.

\begin{remark}

\begin{itemize}
\item Theorem \ref{tp} extends the inequality (\ref{24i-1}) when $F_{N}$ belongs to the Wiener chaos of order $q\geq $ for every $N\geq 1$.  It is known that for every  $N\geq 1$, we have that $\gamma_{N} ^{2} = Var (\frac{1}{q}\Vert  DF_{N} \Vert ^{2} _{H}) $ behaves as $\mathbf{E}F_{N} ^ {4}-3$, more precisely (see Section 5.2.2 in \cite{NPbook}) 
\begin{equation*}
\gamma_{N} ^{2}  \leq \frac{q-1}{3q} (\mathbf{E}F_{N} ^ {4}-3) \leq (q-1) \gamma_{N}^{2}  .
\end{equation*}
It also extends some results in \cite{HLN}  (in particular Theorem 6.2 and Corollary  6. 6). 
 
\item The appearance of truncation function $\Psi$ in the right-hand side of (\ref{p-p}) replaces the condition (\ref{neg})  on the existence of the negative moments of the Malliavin derivative.  
\end{itemize}

\end{remark}

As a consequence of Theorem \ref{tp} we will obtain the asymptotic expansion for the expectation of measurable functionals of $F_{N}$.

\begin{theorem}\label{exp}
Let $p_{N}$ be given by (\ref{pn}) and fix $M\geq 0, \gamma>0$.   Then 
\begin{eqnarray}\label{26i-1}
\sup_{ f\in \mathcal{E}(M, \gamma)} \left| \mathbf{E}f(F_{n})- \int_{\mathbb{R}^{d}} f(x) p_{N} (x)dx \right| = o(\gamma_{N})
\end{eqnarray}
where $\mathcal{E}(M, \gamma)$ is the set of measurable functions $f:\mathbb{R} ^{d}\to \mathbb{R}$ such that $\vert f(x)\vert \leq M (1+\vert x\vert ^{\gamma} )$ for every $x\in \mathbb{R} ^{d}$.

\end{theorem}
{\bf Proof: } For any measurable function $f$ that satisfies the assumptions in the statement, we can write
\begin{eqnarray*}
\mathbf{E}f(F_{n})- \int_{\mathbb{R}^{d}} f(x) p_{N} (x)dx&=& \mathbf{E}f(F_{N})- \mathbf{E}f(F_{N}) \Psi_{N}  +  \mathbf{E} f(F_{N}) \Psi_{N} - \int_{\mathbb{R}^{d}} f(x) p_{N} (x)dx\\
&=&\mathbf{E} f(F_{N})(1-\Psi_{N}  ) + \int_{\mathbb{R}^{d}} f(x) \left[ p^ {\psi} _{F_{N} } (x)- p_{N} (x) \right]dx
\end{eqnarray*}
with $p^ {\psi} _{F_{N} } $ defined in (\ref{pfn}). Then for every $p>1$ and for $  \alpha\geq 0$ large enough (it may depend on $M, \gamma$) with $p'$ such that $\frac{1}{p} + \frac{1}{p' } =1$, we have, by using H\"older's inequality,
\begin{eqnarray*}
&&\left| \mathbf{E}f(F_{n})- \int_{\mathbb{R}^{d}} f(x) p_{N} (x)dx \right|\\
& \leq& \left( \mathbf{E}\left| f(F_{N}) \right| ^ {p' } \right) ^ {\frac{1}{p' }} \Vert 1-\Psi_{N} \Vert _{p}+ \int_{\mathbb{R}^{d}} \vert f(x)\vert \left( \frac{1}{ 1+\vert x\vert ^ {2}}\right) ^ {\frac{{\alpha}}{2}}dx \times  \sup_{x\in \mathbb{R}^{d} } \left( 1+\vert x\vert ^ {2} \right) ^ {\frac{{\alpha}}{2}}\left[ p^ {\psi} _{F_{N} } (x)- p_{N} (x) \right]\\
& \leq& \left(\mathbf{E}\left| f(F_{N}) \right| ^ {p' } \right) ^ {\frac{1}{p' }} \Vert 1-\Psi_{N} \Vert _{p}+ \left( c \Vert 1-\Psi_{N}\Vert _{p} + o(\gamma_{N})\right)   \\
&\leq & C\Vert 1-\Psi_{N} \Vert _{p}  + o(\gamma_{N})=o(\gamma_{N}).
\end{eqnarray*}
where we used (\ref{p-p}). So the desired conclusion is obtained. \qed

We insert some remarks related to the above result.

\begin{remark}
\begin{itemize}
\item 
Theorem \ref{exp} is related to Theorem 9.3.1 in \cite{NPbook} (which reproduces the main result of \cite{NoPe2}).  Notice that we do not assume the differentiability of $f$ and our result includes the uniform control in the case of the Kolmogorov distance.

\item   In particular, the asymptotic expansion (\ref{26i-1}) holds for any bounded function $f$ and for polynomial functions, but the class of examples is obviously larger.

\end{itemize}
\end{remark}

In the case of sequences in a Wiener chaos of fixed order, we have the  following result, which follows from the above results and from the hypercontractivity property (\ref{hyper}).

\begin{corollary}
Assume $(F_{N})_{N\geq 1}=( I_{q}( f_{N}))_{N\geq 1}$ (with  $f_{N} \in H ^{\otimes q}$ for every $N\geq 1$) be a random sequence in the  $q$ th Wiener chaos ($q\geq 2$). Assume $\mathbf{E} F_{N} ^{2}=1 $ for every $N\geq 1$ and suppose that the condition {\bf  [C1\!\! ]} holds with $\gamma_{N}= Var \left( \frac{\Vert DF_{N}\Vert _{H} ^{2}}{q}\right)$ i.e.
$$\left( F_{N}, \gamma_{N}^{-1} \left(  \frac{\Vert DF_{N}\Vert _{H}^{2}}{q} -1\right)  \right)\to ^{(d)}_{N\to \infty}  (Z_{1}, Z_{2})$$
where $(Z_{1}, Z_{2}) $ is a centered Gaussian vector with $\mathbf{E}Z_{1} ^{2}= \mathbf{E}Z_{2} ^{2}=1 $ and  $\mathbf{E}Z_{1} Z_{2}= \rho \in (0,1)$.
Then we have the asymptotic expansions (\ref{p-p}) and (\ref{26i-1}). 
 \end{corollary}
{\bf Proof: } The results follows from Theorems \ref{tp} and \ref{exp} by noticing that (\ref{hyper}) ensures that the conditions {\bf  [C2\!\! ]} and {\bf  [C3\!\! ]} are satisfied. \qed

\section{Examples}
In the last part of our work, we will present three examples in order two illustrate our theory. The first example concerns a one-dimensional sequence of random variables   in the second Wiener chaos and it is inspired from \cite{NPbook}, Chapter 9.  The second example is a multidimensional one and it involves quadratic functions of two  correlated fractional Brownian motions while the last example involves a sequence in a finite sum of Wiener chaoses.

\subsection{A one-dimensional in a Wiener chaos of a fixed order: Quadratic variations of stationary Gaussian sequences}
We consider a centered stationary Gaussian sequence $(X_{k}) _{k \in \mathbb{Z}} $ with  $\mathbf{E} X_{k} ^{2}= 1$ for every $k\in \mathbb{Z} $. We set

$$\rho (l):= \mathbf{E}X_{0} X_{l}  \mbox{ for every } l \in \mathbb{Z}$$
so $\rho(0)=1$ and $\vert \rho(l)\vert \leq 1$ for every $l\in \mathbb{Z}$. 

Without loss of generality and in order to fit to our context, we assume that $X_{k}= W(h_{k})$ with $\Vert h_{k}\Vert _{H}=1$  for every $k\in \mathbb{Z}$, where $(W(h), h\in H) $ is an isonormal process as defined in Section \ref{app}.  Define the quadratic variation sequence
\begin{equation}\label{fnn}
V_{N} =\frac{1}{\sqrt{N}}\sum_{k=0} ^{N}  \left( X_{k} ^{2}-1\right), \hskip0.5cm N\geq 1.
\end{equation}
 and let $v_{N}: = \mathbf{E} V_{N} ^{2}$.  Assume $\rho \in \ell ^{2}(\mathbb{Z})$ (the Banach space of $2$-summable functions endowed with the norm $\Vert \rho\Vert _{\ell ^{2} (Z) }^{2} = \sum _{k\in \mathbb{Z} }\rho( k) ^{2}$). In this case,  we know from \cite{NPbook} that the deterministic sequence $v_{N}$ converges to the constant $2 \sum _{k\in \mathbb{Z} }\rho( k) ^{2}>0$ and consequently it plays no role in the sequel.

Define the renormalized sequence
\begin{equation}
\label{vnn}
F_{N} = \frac{V_{N} }{\sqrt{v_{N}}} = I_{2} \left( \frac{1}{ \sqrt{Nv_{N}}} \sum_{k=0} ^{N} h_{k} ^{\otimes 2}\right) \mbox{ for every } N\geq 1.
\end{equation}
Obviously $ \mathbf{E}F_{N}=0$ and $\mathbf{E}F_{N} ^{2}= 1$ for every $N\geq 1$. It is well known (see \cite{BM}, see also Chapter 7 in \cite{NPbook})  that, as $N \to \infty$ 
\begin{equation*}
F_{N} \to ^{(d)} Z\sim N(0,1). 
\end{equation*}

Let us see that  the sequence $(F_{N}) _{N\geq 1} $ satisfies the conditions  conditions \\{\bf  [C1\!\! ]}-{\bf  [C3\!\! ]} with
\begin{equation}
\label{g1}
\gamma_{N}^{2}= Var \left( \frac{1}{2} \Vert DF_{N} \Vert  _{H} ^{2} \right).
\end{equation}  
Actually, this follows from the computations contained in Section 5 of \cite{NPbook}, but we recall the main ideas because they are needed later. Notice that $\langle DF_{N}, D(-L) ^{-1} F_{N} \rangle _{H} =\frac{1}{2} \Vert DF_{N} \Vert  _{H}^{2}$ since $F_{N}$ belongs to the second Wiener chaos.

To prove {\bf  [C1\!\! ]}, we need to assume 
\begin{equation}
\label{c3}
\rho \in \ell ^{\frac{4}{3} } (\mathbb{Z}).
\end{equation}
 In this case  $\gamma ^{2}   _{N}$ (which coincides with the fourth cumulant of $F_{N}$, denoted $k_{4}(F_{N})$)
 behaves as $c\frac{1}{N}$ for $N$ large, see Theorem 7.3.3 in \cite{NPbook}. The proof of Theorem 9.5.1 in \cite{NPbook} implies that,   as $N\to \infty,$ 
\begin{equation*}
\left( F_{N}, \gamma_{N} ^{-1} \left(  \frac{\Vert DF_{N}\Vert _{H} ^{2}}{2}-1\right)\right)= \left( F_{N}, \gamma_{N} ^{-1} (\langle DF_{N}, D(-L)^{-1} F_{N}\rangle _{H}-1\right)\to  ^{(d) } (Z_{1}, Z_{2} )
\end{equation*}
where $(Z_{1}, Z_{2}) $ is a correlated centered  Gaussian  vector with $\mathbf{E}Z_{1}^{2}= \mathbf{E} Z_{2} ^{2} =1$. Therefore  condition {\bf  [C1\!\! ]} is fullfiled.  

  Conditions {\bf  [C2\!\! ]} and {\bf  [C3\!\! ]} are satisfied since $F_{N}$ belongs to the second Wiener chaos and one can use  the hypercontractivity property of multiple integrals (\ref{hyper}). Consequently  Theorems \ref{exp} and \ref{tp} can be then applied to the sequence (\ref{fnn}).

\begin{remark}
Let $(B^{H}_{t}) _{t\in \mathbb{R}} $ be a (two-sided)  fractional Brownian motion with Hurst parameter $H\in (0,1)$. We recall that $B ^{H}$ is a centered Gaussian process with covariance function 
\begin{equation*}
\mathbf{E} B ^{H}_{t}B ^{H} _{s}= \frac{1}{2} (\vert t\vert ^{2H}+\vert s\vert ^{2H}-\vert t-s\vert ^{2H} )
\end{equation*} 
for every $s,t \in \mathbb{R}$.  Let $X_{k}= B ^{H}_{k+1}- B^{H}_{k}$ for every $k\in \mathbb{Z}$.  In this case $\rho_{H}  (k): =\mathbf{E} X_{0}X_{k} $ is given by 
\begin{equation}
\label{roh}
\rho _{H}(k) = \frac{1}{2} \left( \vert k+1\vert ^{2H} + \vert k-1\vert ^{2H}-2 \vert k\vert ^{2H} \right), \hskip0.5cm \forall k\in \mathbb{Z} .
\end{equation}
Since $\rho _{H}(k)$ behaves as $H(2H-1) \vert k\vert ^{2H-2}$ for $\vert k\vert $ large enough, condition (\ref{c3}) is fullfiled for $0<H<\frac{5}{8}$. Therefore  Theorems \ref{tp} and \ref{exp} can be applied for $H\in (0, \frac{5}{8})$.
\end{remark}

\subsection{A multidimensional: quadratic variations of correlated fractional Brownian motions}

For every $H\in (0,1)$, denote by
\begin{equation*}f^{H} _{t} (u)= d(H) \left( (t-u) _{+} ^{H-\frac{1}{2}}- (-u) _{+} ^{H-\frac{1}{2}}\right), \hskip0.4cm t\in \mathbb{R}
\end{equation*}
with $d(H)$ a positive constant that ensures that $ \int_{\mathbb{R}} f^{H} _{t} (u) ^{2}du = \vert t \vert ^{2} $ for every $t\in \mathbb{R}$.  For every $k\in \mathbb{Z}$, let

\begin{equation}
\label{lk}L _{k} ^{H}= f_{k+1} ^{H}-f_{k+1} ^{H}.
\end{equation}
Consider  $(W_{t}) _{t\in \mathbb{R}}$  a Wiener process on the whole real line and define the two fractional Brownian motions 

\begin{equation*}
B _{t} ^{ H_{i}}=\int_{\mathbb{R} }  f_{t} ^{H_{i}} (s)dW_{s}, \hskip0.5cm i=1,2.
\end{equation*}
The fBms $B^{H_{1}}$ and $ B^{H_{2}}$ are correlated. We actually have (see \cite{MaTu}, Section 4), for every $k,l \in \mathbb{Z}$,

\begin{eqnarray}
&&\mathbf{E} (B ^{H_{1}}_{k+1}- B ^{H_{1}}_{k})(B ^{H_{2}}_{l+1}- B ^{H_{2}}_{l})= \langle f^{H_{1}}_{k+1}-  f^{H_{1}}_{k},  f^{H_{2}}_{l+1}-f^{H_{2}}_{l}\rangle =\langle L ^{H_{1}}_{k}, L ^{H_{2}}_{l} \rangle \nonumber \\
&=& D(H_{1}, H_{2}) \rho _{ \frac{H_{1}+ H_{2}}{2}}(k-l) \label{ex-1}
\end{eqnarray}
where $ D(H_{1}, H_{2})$ is an explicit constant such that $D(H_{i}, H_{i})=1$ for $i=1,2$.

Assume $H_{1}, H_{2} \in \left(0, \frac{3}{4}\right). $ Define the quadratic variations, for $i=1,2$ and for $N\geq1$

\begin{equation*}
V_{N} ^{(i) } = \frac{1}{\sqrt{N}}\sum _{k=0} ^{N-1}\left[ \left( B_{k+1} ^{H_{i} }-B_{k} ^{H_{i} }\right) ^{2} -1\right]
\end{equation*} 
and
\begin{equation}\label{24m-1}
F_{N} ^{(i)} =\frac{V_{N} ^{(i)}}{ \sqrt{v_{n} ^{(i) }}}
\end{equation}
with $v_{N} ^{(i) }= \mathbf{E} (V_{N} ^{(i)}) ^{2}$, $i=1,2$. Recall that $v_{N} ^{(i)}\to _{N\to \infty} c_{v_{i}}:= 2\sum_{k\in \mathbb{Z}} \rho_{H_{i}} (k)^{2}$ for $i=1,2$.

Let us show that the sequence $F_{N}= (F_{N} ^{(1)}, F_{N} ^{(2)})$ satisfies the conditions\\ {\bf  [C1\!\! ]}-{\bf  [C3\!\! ]} with $\gamma_{N} = N ^{-\frac{1}{2}}$. We  will focus on {\bf  [C1\!\! ]} (the assumptions {\bf  [C2\!\! ]} and {\bf  [C3\!\! ]} can be easily checked since we deal with elements of the second Wiener chaos and we can apply the hypercontractivity (\ref{hyper})).

From the previous paragraph we know that for $i=1,2$, the random vectors

\begin{equation}
\label{25m-1}
\left( F_{N} ^{(i)}, \sqrt{N}\langle DF_{N} ^{(i)}, D(-L) ^{-1} F_{N} ^{(i)}\rangle_{H}-1\right) 
\end{equation}
converge in distribution, when $N\to \infty$, to some  centered  Gaussian vectors $(Z_{i} ^{(1)}, Z_{i} ^{(2)})$ with $\mathbf{E} (Z_{i} ^{(1) }  ) ^{2}=C_{i,i}$, $i=1,2$ and with nontrivial correlation.

The main tool to prove the the two-dimensional sequence $(F_{N}^{ (1)}, F_{N}^{(2)} )$ satisfies {\bf  [C1\!\! ]} is the multidimensional version of the Fourth Moment Theorem. Let us recall this result which says that for sequences of random vectors on Wiener chaos, the componentwise convergence to the Gaussian law implies the joint convergence (see \cite{PeTu} or Theorem 6.2.3 in \cite{NPbook}).
\begin{theorem}\label{m4m}
Let $d\geq 2$ and $q_{1}, ..., q_{d} \geq 1$ integers. Consider the random vector
$$F_{N}= \left( F_{N} ^{(1)},\ldots , F_{N} ^{(d)}\right)= \left( I_{q_{1}}(f_{N} ^{(1)}), \ldots , I_{q_{d}}(f_{N} ^{(d)})\right)$$
with $ f_{N} ^{(i)}\in H ^{\otimes q_{i}}$, $i=1,..,d$. Assume that
\begin{equation}
\label{25m-4}
\mathbf{E} F_{N} ^{(i)} F_{N} ^{(j)}\to _{N\to \infty} C(i,j) \mbox{ for every }i, j=1,..,d.
\end{equation} Then the following two conditions are equivalent:
\begin{itemize}
\item $F_{N}$ converges in law, as $N\to \infty$,  to a $d$ dimensional centered Gaussian vector with covariance matrix $C.$
\item For every $i=1,..,d$, $F_{N} ^{(i)}$ converges in law, as $N\to \infty$,  to $N(0, C(i,i))$.
\end{itemize}
\end{theorem}

We will need the two auxiliary lemmas below. The first concerns the asymptotic covariance of $F_{N} ^{(1)} $ and  $F_{N}^{(2)}$. 

\begin{lemma}\label{lex-2}
Let $F_{N} ^{(1) }, F_{N} ^{(2)} $ be defined by (\ref{24m-1}). Assume $H_{1}, H_{2}\in \left( 0, \frac{3}{4}\right).$ Then for every $H_{1}, H_{2} \in (0,1)$, 
\begin{equation*}
\mathbf{E}F_{N} ^{(1)} F_{N} ^{(2)}\to _{N\to \infty} C_{1,2}
\end{equation*}
with
\begin{equation}
\label{c12}
C_{1,2}=  \frac{2 D ^{2} (H_{1}, H_{2})}{\sqrt{c_{v_{1}}c_{v_{2}}}}\sum_{l\in \mathbb{Z} } \rho_{\frac{H_{1}+H_{2}}{2}}(l) ^{2}.
\end{equation}
Consequently, 
\begin{equation}\label{1n-4}
\left( F_{N} ^{(1)}, F_{N} ^{(2)}\right) \to ^{(d)}_{N\to \infty} (Z_{1}^{(1)}, Z_{1}^{(2)}) 
\end{equation} 
where $(Z_{1}^{(1)}, Z_{1}^{(2)})$ is a Gaussian vector with $ \mathbf{E} Z_{1} ^{(1)} =\mathbf{E} Z_{1} ^{(2)}=0$ and 
$$\mathbf{E} ( Z_{1} ^{(1)} ) ^{2}= C_{1,1}=1, \hskip0.2cm \mathbf{E} ( Z_{1} ^{(2)} ) ^{2} =C_{2,2}>0,\hskip0.2cm \mathbf{E}  ( Z_{1} ^{(1)} Z_{1} ^{(2)} )=C_{1,2}$$
with $C_{1,2}$ given by (\ref{c12}). Recall that we denoted by $" \to ^{(d)}" $ the convergence in distribution.

\end{lemma}
{\bf Proof: } By the chaos expression of $F_{N} ^{(i)}$  (\ref{vnn}) and by using the isometry (\ref{iso})
\begin{eqnarray*}
\mathbf{E}F_{N} ^{(1)} F_{N} ^{(2)}&=& \frac{1}{N \sqrt{v_{N} ^{(1)} v_{N} ^{(2)}}} \mathbf{E} I_{2} \left( \sum_{i=0} ^{N-1} (L_{i} ^{H_{1}}) ^{\otimes 2} \right)I_{2}  \left(\sum_{j=0} ^{N-1} (L_{j} ^{H_{2}}) ^{\otimes 2}\right)\\
&=& \frac{1}{N \sqrt{v_{N} ^{(1)} v_{N} ^{(2)}}} 2 \sum_{i,j=0} ^{N-1} \langle L_{i} ^{H_{1}}, L _{j} ^{H_{2}}\rangle  ^{2} 
\end{eqnarray*}
with $L^{H}$ from (\ref{lk}).  From the identity (\ref{ex-1})
\begin{eqnarray*}
\mathbf{E}F_{N} ^{(1)} F_{N} ^{(2)}
&=&\frac{1}{N \sqrt{v_{N} ^{(1)} v_{N} ^{(2)}}} 2 D ^{2} (H_{1}, H_{2})  \sum_{i,j=0} ^{N-1}\rho_{\frac{H_{1}+H_{2}}{2}} (i-j) ^{2} \\
&=&  \frac{2}{ \sqrt{v_{N} ^{(1)} v_{N} ^{(2)}}} \sum_{l=0} ^{N-1} \rho_{\frac{H_{1}+H_{2}}{2}}(l) ^{2} (1-\frac{l}{N}) \\
&\to _{N\to \infty}& \frac{2 D ^{2} (H_{1}, H_{2})}{\sqrt{c_{v_{1}}c_{v_{2}}}}\sum_{l\in \mathbb{Z} } \rho_{\frac{H_{1}+H_{2}}{2}}(l) ^{2} .
\end{eqnarray*}
The sum $\sum_{l\in \mathbb{Z} } \rho_{\frac{H_{1}+H_{2}}{2}}(l) ^{2} $ is finite since $H_{1}+ H_{2}<\frac{3}{2}$. The limit (\ref{1n-4}) will then follow from Theorem \ref{m4m}. \qed

\vskip0.2cm 
We will also need to control the speed of convergence of the covariance function of $F_{N}^{(1)}$ and $F_{N} ^{(2)}$.  In this case we need to impose an additional restriction on the Hurst parameters.

\begin{lemma}\label{f12}
Let $F_{N} ^{(1) }, F_{N} ^{(2)} $ be defined by (\ref{24m-1}) and let $C_{1,2}$ be given by (\ref{c12}). Assume $H_{1}, H_{2}\in \left( 0, \frac{5}{8}\right). $
Then 
$$\sqrt{N} ( \mathbf{E}F_{N} ^{(1)} F_{N} ^{(2)}- C_{1,2} ) \to _{N\to \infty} 0.$$

\end{lemma}
{\bf Proof: } For every $N\geq 1$ we can write 
\begin{eqnarray*}
&&\sqrt{N} ( \mathbf{E}F_{N} ^{(1)} F_{N} ^{(2)}- C_{1,2} )\\
&=&
2D^{2} (H_{1}, H_{2}) \sqrt{N}\left[ \frac{1}{ \sqrt{v_{N} ^{(1)} v_{N} ^{(2)}}} \sum_{l=0} ^{N-1} \rho_{\frac{H_{1}+H_{2}}{2}}(l) ^{2} (1-\frac{l}{N}) -\frac{1}{\sqrt{c_{v_{1}}c_{v_{2}}}}\sum_{l\in \mathbb{Z} } \rho_{\frac{H_{1}+H_{2}}{2}}(l) ^{2} \right].
\end{eqnarray*}
Thus, for $N$ large enough, using that $v_{N} ^{(i)}$ converges to $c_{v_{i}}$ for $i=1,2$, we can write
\begin{eqnarray*}
&&\sqrt{N} \left| \mathbf{E}F_{N} ^{(1)} F_{N} ^{(2)}- C_{1,2} \right| \\
&\leq & C \sqrt{N} \sum _{l \geq N}\rho_{\frac{H_{1}+H_{2}}{2}}(l) ^{2}  +C\frac{1}{\sqrt{N}} \sum_{l=0} ^{N-1} l\times \rho_{\frac{H_{1}+H_{2}}{2}}(l) ^{2} + C \sqrt{N}\left|\frac{1}{ \sqrt{v_{N} ^{(1)} v_{N} ^{(2)}}}- \frac{1}{\sqrt{c_{v_{1}}c_{v_{2}}}}\right|\\
&:= &T_{1,N}+T_{2,N}+ T_{3,N}.
\end{eqnarray*}
Since for $l$ large, the function $\rho_{H} (l) $ behaves as $C \times \vert l\vert ^{2H-2}$, we get, since $H_{1}, H_{2} <\frac{5}{8}$ (so $H_{1}+ H_{2}< \frac{5}{4}$)
\begin{equation*}
T_{1,N}\leq C \sqrt{N} N ^{2H_{1}+ 2H_{2}-3} =C N ^{ 2H_{1}+ 2H_{2}- \frac{5}{2}}\to _{N\to \infty } 0
\end{equation*}
and 
\begin{equation*}
T_{2,N}\leq C \frac{1}{\sqrt{N}} N ^{2H_{1}+ 2H_{2}-2} =C N ^{ 2H_{1}+ 2H_{2}- \frac{5}{2}}\to _{N\to \infty } 0.
\end{equation*}
Concerning the summand $T_{3,N}$,  we have for $N$ large enough
\begin{eqnarray*}
T_{3,N} &\leq & C \sqrt{N} \left| \frac{1}{\sqrt{v_{N} ^{(1)}}}- \frac{1}{\sqrt{c_{v_{1}}}}\right|+C \sqrt{N} \left| \frac{1}{\sqrt{v_{N} ^{(2)}}}- \frac{1}{\sqrt{c_{v_{2}}}}\right|\\
&\leq &
C\sqrt{N}\left( \vert v_{N}^{(1)}- c_{v_{1}}\vert + \vert v_{N}^{(2)}- c_{v_{2}}\vert \right) \\
&\leq &C\sqrt{N}\left(   \sum_{l\geq N} \rho ^{2}_{H_{1}} (l) + \sum_{l\geq N} \rho ^{2}_{H_{2}} (l)\right)\leq  C (N ^{4H_{1}-\frac{5}{2}}+ N ^{4H_{2}-\frac{5}{2}})\to_{N\to \infty }0
\end{eqnarray*}
and the conclusion follows. \qed

\vskip0.2cm

The following lemma is also important in order to check {\bf  [C1\!\! ]}.

\begin{lemma}\label{ex-3}Assume $H_{1}<\frac{5}{8}$ and $H_{2}<\frac{5}{8}$. Let $F_{N} ^{(1) }, F_{N} ^{(2)} $ be defined by (\ref{24m-1}) and let $C_{1,2}$ be given by (\ref{c12}). Then
\begin{equation*}\sqrt{N}\left( \frac{1}{2} \langle DF_{N} ^{(1)}, DF_{N} ^{(2)}\rangle_{H}  - C_{1,2}\right) 
\end{equation*}
converges in distribution, as $N$ goes to infinity, to a Gaussian random variable.
\end{lemma}
{\bf Proof: } Since $\mathbf{E} F_{N}^{(1)} F_{N} ^{(2)} = \frac{1}{2} \langle DF_{N} ^{(1) }, DF_{N} ^{(2) }\rangle _{H}$ and by Lemma \ref{f12}
$$\sqrt{N} \left( \mathbf{E} \frac{1}{2}\langle DF_{N} ^{(1)}, DF_{N} ^{(2)}\rangle _{H} -C_{1,2}\right) \to _{N\to \infty}0$$
it suffices to show that the sequence $(Y_{N}) _{N\geq 1} $ with 

\begin{equation}\label{yn}
Y_{N}: = \frac{1}{2} \sqrt{N} \left( \langle DF_{N} ^{(1)}, DF_{N} ^{(2)}\rangle _{H} -\mathbf{E}\langle DF_{N} ^{(1)}, DF_{N} ^{(2)}\rangle _{H}\right) 
\end{equation}
converges in distribution as $N$ goes to infinity, to a Gaussian random variable. We will apply the Fourth Moment Theorem. In order to do it, we need to show that 
\begin{equation}
\label{24m-4}
\mathbf{E} Y_{N} ^{2}\to _{N\to \infty } C_{0}>0 \mbox{ and } Var (\Vert DY_{N} \Vert _{H}^{2} ) \to_{N\to \infty } 0.
\end{equation}
From (\ref{vnn}) and the product formula for multiple integrals (\ref{prod}), we can write
\begin{eqnarray}
Y_{N} &=& \frac{1}{\sqrt{N v_{N} ^{(1)} v_{N} ^{(2)}}} I_{2} \left( \sum_{k, l =0} ^{N-1} (L_{k} ^{H_{1}}\otimes L_{l} ^{H_{2}}) \langle L_{k} ^{H_{1}},  L_{l} ^{H_{2}}\rangle \right) \nonumber\\
&=& D(H_{1}, H_{2}) \frac{1}{\sqrt{N v_{N} ^{(1)} v_{N} ^{(2)}}} I_{2} \left( \sum_{k, l =0} ^{N-1} \rho_{\frac{H_{1}+H_{2}}{2}}(k-l) (L_{k} ^{H_{1}}\otimes L_{l} ^{H_{2}}) \right) \label{5m-1}
\end{eqnarray}
with $L^{H} _{k}$ defined by (\ref{lk}).  Consequently, since any four functions of two variable $f,g,u,v$ we have

$$\langle f\tilde{\otimes }g, u\tilde{\otimes } v\rangle = \frac{1}{2}( \langle f,u\rangle \langle g,v\rangle + \langle f,v\rangle \langle g,u\rangle)$$
we will obtain 

\begin{eqnarray}
\mathbf{E}Y_{N} ^{2}&=& D^{2}(H_{1}, H_{2} ) \frac{1}{Nv_{N} ^{(1)}v_{N} ^{(2)} }\sum_{k,l,i,j=0}^{N-1} \rho_{\frac{H_{1}+H_{2}}{2}}(k-l) \rho_{\frac{H_{1}+H_{2}}{2}}(i-j) \langle L_{k} ^{H_{1}}\tilde{\otimes } L_{l} ^{H_{2}}, L_{i} ^{H_{1}}\tilde{\otimes } L_{j} ^{H_{2}}\rangle \nonumber\\
&=& D^{2}(H_{1}, H_{2} ) \frac{1}{Nv_{N} ^{(1)}v_{N} ^{(2)} }\sum_{k,l,i,j=0}^{N-1} \rho_{\frac{H_{1}+H_{2}}{2}}(k-l) \rho_{\frac{H_{1}+H_{2}}{2}}(i-j) \nonumber\\
&&\times \frac{1}{2} \left[ \rho_{H_{1}} (k-i) \rho_{H_{2}}(l-j) + D(H_{1}, H_{2}) ^{2}\rho_{\frac{H_{1}+H_{2}}{2}}(i-l) \rho_{\frac{H_{1}+H_{2}}{2}}(k-j) \right].\label{24m-3}
\end{eqnarray}
For two sequences $(u(n), n\in \mathbb{Z})$ and $(v(n), n\in\mathbb{Z})$, we define their convolution by
$$(u\ast v)(j)= \sum_{n\in \mathbb{Z}}u(n) v(j-n).$$ We will need the Young 's inequality: if $s,p,q\geq 1$ with $\frac{1}{s}+1= \frac{1}{p}+\frac{1}{q}$ then
\begin{equation}
\label{24m-2}
\Vert u\ast v \Vert _{\ell ^{s} (\mathbb{Z})}\leq \Vert u\Vert _{\ell ^{p} (\mathbb{Z})}\Vert v\Vert _{\ell ^{q} (\mathbb{Z})}.
\end{equation}
By the proof of Theorem 7.3.3 in \cite{NPbook}, we have

$$\frac{1}{N} \sum_{k,l,i,j=}^{N-1} \rho_{\frac{H_{1}+H_{2}}{2}}(k-l) \rho_{\frac{H_{1}+H_{2}}{2}}(i-j)\rho_{\frac{H_{1}+H_{2}}{2}}(i-l) \rho_{\frac{H_{1}+H_{2}}{2}}(k-j) \to _{N\to \infty} \langle \rho_{\frac{H_{1}+H_{2}}{2}} ^{\ast 3}, \rho_{\frac{H_{1}+H_{2}}{2}}\rangle _{\ell^{2} (\mathbb{Z})} $$
where $\langle \rho_{\frac{H_{1}+H_{2}}{2}} ^{\ast 3}, \rho_{\frac{H_{1}+H_{2}}{2}}\rangle _{\ell^{2} (\mathbb{Z})}<\infty$ when $H_{1}+ H_{2} \leq \frac{5}{4}$ which implies  that $ \rho_{\frac{H_{1}+H_{2}}{2}}\in \ell ^{\frac{4}{3}}(\mathbb{Z}).$ The second summand in (\ref{24m-3}) can be handled as follows
\begin{eqnarray}
&&\frac{1}{N} \sum_{k_{1},k_{2}, k_{3}, k_{4} =0} ^{N-1} \rho_{\frac{H_{1}+H_{2}}{2}}(k_{4}-k_{3})\rho_{\frac{H_{1}+H_{2}}{2}}(k_{2}-k_{1})\rho_{H_{1}}(k_{1}-k_{4})\rho_{H_{2}} (k_{3}-k_{2})\nonumber \\
&=&\frac{1}{N}\sum_{k_{1}=0} ^{N-1} \sum_{k_{2}, k_{3}, k_{4} =-k_{1}} ^{N-1-k_{1}} \rho_{\frac{H_{1}+H_{2}}{2}}(k_{4}-k_{3})\rho_{\frac{H_{1}+H_{2}}{2}}(k_{2})\rho_{H_{1}}(k_{4})\rho_{H_{2}} (k_{3}-k_{2})\nonumber \\
&=&\frac{1}{N}\sum_{k_{2}, k_{3}, k_{4} \in \mathbb{Z}}\rho_{\frac{H_{1}+H_{2}}{2}}(k_{4}-k_{3})\rho_{\frac{H_{1}+H_{2}}{2}}(k_{2})\rho_{H_{1}}(k_{4})\rho_{H_{2}} (k_{3}-k_{2})\nonumber\\
&&\left[ 1\vee \left(1-\frac{\max(k_{2}, k_{3}, k_{4})}{N}\right) -0\wedge \left( \frac{\min (k_{2}, k_{3}, k_{4})}{N}\right)\right]1_{\vert k_{2} \vert <N, \vert k_{3} \vert <N, \vert k_{4} \vert <N}\nonumber \\
&&\to _{N\to \infty }\sum_{k_{2}, k_{3}, k_{4} \in \mathbb{Z}}\rho_{\frac{H_{1}+H_{2}}{2}}(k_{4}-k_{3})\rho_{\frac{H_{1}+H_{2}}{2}}(k_{2})\rho_{H_{1}}(k_{4})\rho_{H_{2}} (k_{3}-k_{2})\nonumber \\
&=& \langle \rho_{\frac{H_{1}+H_{2}}{2}}\ast \rho_{H_{1}}\ast \rho_{H_{2}}, \rho_{\frac{H_{1}+H_{2}}{2}}\rangle _{\ell ^{2} (\mathbb{Z})}.\label{25m-3}
\end{eqnarray}
Again notice that  $\langle \vert \rho_{\frac{H_{1}+H_{2}}{2}}\vert\ast \vert  \rho_{H_{1}}\vert \ast  \vert  \rho_{H_{2}}\vert, \vert \rho_{\frac{H_{1}+H_{2}}{2}}\vert \rangle _{\ell ^{2} (\mathbb{Z})}<\infty$  since by H\"older's  and Young's inequalities, by using  $\rho_{\frac{H_{1}+H_{2}}{2}}\in \ell ^{\frac{4}{3}}(\mathbb{Z})$,

\begin{eqnarray*}
\langle \vert \rho_{\frac{H_{1}+H_{2}}{2}}\vert\ast \vert  \rho_{H_{1}}\vert \ast  \vert  \rho_{H_{2}}\vert, \vert \rho_{\frac{H_{1}+H_{2}}{2}}\vert \rangle _{\ell ^{2} (\mathbb{Z})}
&\leq& C \Vert \vert \rho_{\frac{H_{1}+H_{2}}{2}}\vert\ast \vert  \rho_{H_{1}}\vert \ast  \vert  \rho_{H_{2}}\vert \Vert  _{\ell ^{4} (\mathbb{Z})}\\
&\leq & C \Vert \vert  \rho_{\frac{H_{1}+H_{2}}{2}} \vert \Vert _{\ell ^{\frac{4}{3}}(\mathbb{Z})}  \times \Vert \vert  \rho_{H_{1}}\vert \ast  \vert  \rho_{H_{2}}\vert \Vert _{\ell ^{2} (\mathbb{Z})}\\
&\leq&  \Vert  \rho_{\frac{H_{1}+H_{2}}{2}} \Vert _{\ell ^{\frac{4}{3}}(\mathbb{Z})} \times  \Vert \vert  \rho_{H_{1}}\vert  \Vert _{\ell ^{\frac{4}{3}}(\mathbb{Z})} \times  \Vert \vert  \rho_{H_{2}}\vert  \Vert _{\ell ^{\frac{4}{3}}(\mathbb{Z})}
\end{eqnarray*}
where we applied twice (\ref{24m-2}). Thus
$$\mathbf{E}Y_{N}^{2} \to_{N\to \infty } C_{0}:= \frac{D^{2} (H_{1}, H_{2}) }{c _{v_{1}}c_{v_{2}}}\frac{1}{2} \left(  \langle \rho_{\frac{H_{1}+H_{2}}{2}}\ast \rho_{H_{1}}\ast \rho_{H_{2}}, \rho_{\frac{H_{1}+H_{2}}{2}}\rangle _{\ell ^{2} (\mathbb{Z})}+ D^{2} (H_{1}, H_{2}) \langle \rho_{\frac{H_{1}+H_{2}}{2}} ^{\ast 3}, \rho_{\frac{H_{1}+H_{2}}{2}}\rangle _{\ell^{2} (\mathbb{Z})} \right)$$
and the first part of (\ref{24m-4}) is obtained. Next, we check the convergence of the Malliavin derivative in (\ref{24m-4}).  From (\ref{5m-1}) and the product formula  (\ref{prod})

\begin{equation*}
\Vert DY _{N} \Vert _{H}^{2} = \frac{4D (H_{1}, H_{2}) ^{2} }{N v_{N} ^{(1)} v_{N} ^{(2)}} I_{2} \left(   \sum_{i,j,k.l=0}^{N-1} \rho_{\frac{H_{1}+H_{2}}{2}}(k-l)\rho_{\frac{H_{1}+H_{2}}{2}}(i-j)
(L_{k} ^{H_{1}}\tilde{\otimes} L _{l} ^{H_{2}})\otimes _{1}  (L_{i} ^{H_{1}}\tilde{\otimes} L _{j} ^{H_{2}})\right)
\end{equation*}
where the contraction $\otimes _{1} $ is defined in (\ref{contra}). Since for every $f,g, u,v \in L ^{2} (\mathbb{R} ^{2})$, we have 
$$(f\tilde{\otimes } g)\otimes _{1} (u\tilde{\otimes v})= \frac{1}{4} \left( \langle f,u\rangle (g\otimes v) + \langle f,v \rangle (g\otimes u)+ \langle g,u\rangle ( f\otimes v)+ \langle g,v\rangle (f\otimes u)\right) $$
we obtain, via (\ref{ex-1})
\begin{eqnarray*}
\Vert DY _{N} \Vert_{H} ^{2} &=&\frac{D (H_{1}, H_{2}) ^{2} }{N v_{N} ^{(1)} v_{N} ^{(2)}} I_{2} \left(   \sum_{i,j,k.l=0}^{N-1} \rho_{\frac{H_{1}+H_{2}}{2}}(k-l)\rho_{\frac{H_{1}+H_{2}}{2}}(i-j)     \right. \\
&&\left.   \left[ \rho_{H_{1}} (k-i) (L _{l} ^{H_{2}}\otimes L _{j} ^{H_{2}} )+ D(H_{1}, H_{2})\rho_{\frac{H_{1}+H_{2}}{2}}(k-j) (L _{l} ^{H_{2}}\otimes L_{i} ^{H_{1}})\right. \right. \\
&&\left. \left. + D(H_{1}, H_{2}) \rho_{\frac{H_{1}+H_{2}}{2}} (l-i)( L _{k} ^{H_{1}}\otimes L_{j} ^{H_{2}})+ \rho_{H_{2}} (l-j) (L_{k} ^{H_{1}}\otimes L _{i} ^{H_{1}})\right] \right)
\end{eqnarray*}
and by the isometry formula (\ref{iso})

\begin{eqnarray*}
&&Var (\Vert DY_{N} \Vert _{H} ^{2}) \\
&=& C _{1} \frac{1}{N^{2} (v_{N} ^{(1) }) ^{2} (v_{N} ^{(2) }) ^{2} }\sum_{i,j,k.l, i', j', k', l'=0}^{N-1} \rho_{\frac{H_{1}+H_{2}}{2}}(k-l)\rho_{\frac{H_{1}+H_{2}}{2}}(i-j) \rho_{\frac{H_{1}+H_{2}}{2}}(k'-l')\rho_{\frac{H_{1}+H_{2}}{2}}(i'-j')\\ 
&& \times \rho_{\frac{H_{1}+H_{2}}{2}}(k-i) \rho_{\frac{H_{1}+H_{2}}{2}}(k'-i') \rho_{\frac{H_{1}+H_{2}}{2}}(l- l') \rho_{\frac{H_{1}+H_{2}}{2}} (j-j')\\
&&+ C _{2}\frac{1}{N^{2} (v_{N} ^{(1) }) ^{2} (v_{N} ^{(2) }) ^{2} }\sum_{i,j,k.l, i', j', k', l'=0}^{N-1} \rho_{\frac{H_{1}+H_{2}}{2}}(k-l)\rho_{\frac{H_{1}+H_{2}}{2}}(i-j) \rho_{\frac{H_{1}+H_{2}}{2}}(k'-l')\rho_{\frac{H_{1}+H_{2}}{2}}(i'-j')\\ 
&& \times \rho_{H_{1}} (k-i) \rho_{H_{1}}  (k'-i') \rho_{H_{2}} (l- l') \rho_{H_{2}} (j-j ') \\
&&+C _{2}\frac{1}{N^{2} (v_{N} ^{(1) }) ^{2} (v_{N} ^{(2) }) ^{2} }\sum_{i,j,k.l, i', j', k', l'=0}^{N-1} \rho_{\frac{H_{1}+H_{2}}{2}}(k-l)\rho_{\frac{H_{1}+H_{2}}{2}}(i-j) \rho_{\frac{H_{1}+H_{2}}{2}}(k'-l')\rho_{\frac{H_{1}+H_{2}}{2}}(i'-j')\\ 
&& \times \rho_{H_{1}} (k-i) \rho_{\frac{H_{1}+H_{2}}{2}}(k'- j') \rho_{H_{2}} (l- l') \rho_{\frac{H_{1}+H_{2}}{2}}(j-i').
\end{eqnarray*}
We can write, using the convolution symbol
\begin{eqnarray}
&&Var (\Vert DY_{N} \Vert _{H}^{2}) \nonumber \\
&\leq &  C \frac{1}{  N ^{2}} \sum_{l, j'=0}^{N-1} \left[ \left(  \vert  \rho_{\frac{H_{1}+H_{2}}{2}}\vert \ast  \vert  \rho_{\frac{H_{1}+H_{2}}{2}}\vert \ast  \vert  \rho_{\frac{H_{1}+H_{2}}{2}}\vert  \ast  \vert  \rho_{\frac{H_{1}+H_{2}}{2}}\vert \right) (l- j') ^{2}  \right.  \nonumber\\
&&\left. +\left(  \vert  \rho_{\frac{H_{1}+H_{2}}{2}}\vert  \vert  \rho_{\frac{H_{1}+H_{2}}{2}}\vert \ast \vert  \rho_{H_{1}} \vert \ast \vert  \rho_{H_{2}} \vert \right) (l- j') ^{2} \right. \nonumber \\
&&\left. + \left(  \vert  \rho_{\frac{H_{1}+H_{2}}{2}}\vert \ast  \vert  \rho_{\frac{H_{1}+H_{2}}{2}}\vert \ast  \vert  \rho_{\frac{H_{1}+H_{2}}{2}}\vert  \ast  \vert  \rho_{\frac{H_{1}+H_{2}}{2}}\vert \right) (l- j') 
\left(  \vert  \rho_{\frac{H_{1}+H_{2}}{2}}\vert  \vert  \rho_{\frac{H_{1}+H_{2}}{2}}\vert \ast \vert  \rho_{H_{1}} \vert \ast \vert  \rho_{H_{2}} \vert \right) (l- j') \right] \nonumber\\
&\leq & C \frac{1}{N} \sum_{l=-(N-1) } ^{N-1} g(l) \label{11n-1}
\end{eqnarray}
with
$$ g(l) =\left(  \vert  \rho_{\frac{H_{1}+H_{2}}{2}}\vert \ast  \vert  \rho_{\frac{H_{1}+H_{2}}{2}}\vert \ast  \vert  \rho_{\frac{H_{1}+H_{2}}{2}}\vert  \ast  \vert  \rho_{\frac{H_{1}+H_{2}}{2}}\vert \right) (l) ^{2} +\left(  \vert  \rho_{\frac{H_{1}+H_{2}}{2}}\vert \ast   \vert  \rho_{\frac{H_{1}+H_{2}}{2}}\vert \ast \vert  \rho_{H_{1}} \vert \ast \vert  \rho_{H_{2}} \vert \right) (l) ^{2} $$

$$+\left(  \vert  \rho_{\frac{H_{1}+H_{2}}{2}}\vert \ast  \vert  \rho_{\frac{H_{1}+H_{2}}{2}}\vert \ast  \vert  \rho_{\frac{H_{1}+H_{2}}{2}}\vert  \ast  \vert  \rho_{\frac{H_{1}+H_{2}}{2}}\vert \right) (l) 
\left(  \vert  \rho_{\frac{H_{1}+H_{2}}{2}}\vert  \vert  \rho_{\frac{H_{1}+H_{2}}{2}}\vert \ast \vert  \rho_{H_{1}} \vert \ast \vert  \rho_{H_2} \vert \right) (l) .$$

Since for any two sequence $u,v \in \ell ^{2} (\mathbb{Z}) $ we have 
$$\limsup_{n\to \infty}(\vert u \vert \ast \vert v \vert )(n) \leq \sqrt{ \sum _{j\geq M} u(j) ^{2} }\sqrt{ \sum _{j\geq M} v(j) ^{2} }
$$ for every $M\geq 1$ (see \cite{NPbook}, page 168), we obtain as in \cite{NPbook}, by using again  Young's inequality

$$g(N)\to _{N\to \infty }0.$$
The convergence of  $Var (\Vert DY_{N} \Vert _{H}^{2})$ to zero follows from (\ref{11n-1}) and  C\'esaro theorem. \qed

\vskip0.2cm

Let us go back to the sequence $F_{N}= (F_{N} ^{(1)}, F_{N} ^{(2)})$ from  (\ref{24m-1}).
\begin{prop}
The two-dimensional random sequence $F_{N}=\left( F_{N} ^{(1)}, F_{N} ^{(2)} \right)$ defined by (\ref{24m-1}) satisfies condition  {\bf  [C1\!\! ]} with $\gamma_{N}= N ^{-\frac{1}{2}}.$
\end{prop}
{\bf Proof: } Since for every $N\geq 1$
\begin{equation}
\label{11m-2}\langle DF_{N} ^{(1) }, D(-L) ^{-1} F_{N} ^{(2)} \rangle= \langle DF_{N} ^{(2) }, D(-L) ^{-1} F_{N} ^{(1)} \rangle= \frac{1}{2} \langle DF_{N} ^{(1) }, DF_{N} ^{(2) }\rangle
\end{equation}
we need to show that the random vector

\begin{equation}
\label{25m-2}
\left( F_{N} ^{(1)}, F_{N} ^{(2) }, \sqrt{N}\left(\frac{1}{2} \Vert D F_{N} ^{(1) }\Vert_{H} ^{2} -1\right) , \sqrt{N}\left(\frac{1}{2} \Vert D F_{N} ^{(2) }\Vert_{H} ^{2} -1\right) , \sqrt{N} (\frac{1}{2}\langle DF_{N} ^{(1)}, DF_{N} ^{(2)}\rangle_{H}  - C_{1,2})  \right) 
\end{equation}
converges to a Gaussian vector.

All the components of the above vector  are multiple integrals of order 2 or sum of a multiple integral of order 2 and a deterministic term tending to zero (the case of the last component). We will use the multidimensional Fourth Moment Theorem to  prove the convergence of (\ref{25m-2}). From the convergence of the vectors (\ref{25m-1}),  by using  Lemmas \ref{lex-2} and \ref{ex-3}, it suffices to check the condition  (\ref{25m-4})  from the multidimensional Fourth Moment Theoren (Theorem \ref{m4m}), i.e.  the covariances of the components of the vector (\ref{25m-2}) converge to constants as $N\to \infty$. The covariance of $F_{N} ^{(1)}$ and $F_{N} ^{(2) }$ converges to $C_{1,2}$ (\ref{c12}) due to Lemma \ref{lex-2}. Concerning the other limits, we have from (\ref{vnn}) (recall the the symbol $\sim$ means that the sides have the same limit as $N\to \infty$)

\begin{eqnarray*}
\mathbf{E} F_{N } ^ {(1)}\left( \frac{ \Vert DF_{N} ^{(2)}\Vert _{H} ^{2}}{2}-1\right)&\sim& C \mathbf{E} \frac{1}{N}   I_{2}\left( \sum_{i=0} ^{N-1}(L_{k} ^{H_{1}} ) ^{\otimes 2}   \right)I_{2} \left( \sum_{j,k=0} ^{N-1} \rho_{\frac{H_{2}+H_{2}}{2}} (j-k) (L ^{H_{1}}_{j} \otimes L ^{H_{2}}_{k} )\right)\\
&=&C \frac{1}{N} \sum _{i,j,k=0}^{N-1} \rho_{\frac{H_{1}+H_{2}}{2}} (j-k)\rho_{H_{1}} (i-j) \rho_{\frac{H_{1}+H_{2}}{2}} (i-k)\\
&& \to _{N\to \infty } \langle \rho_{\frac{H_{1}+H_{2}}{2}} \ast \rho_{H_{1}}, \rho_{\frac{H_{1}+H_{2}}{2}} \rangle _{\ell ^{2} (\mathbb{Z})}<\infty
\end{eqnarray*}
where the limit is obtained similarly as (\ref{25m-3}) and the last series is finite by using Young inequality and   $\rho_{\frac{H_{1}+H_{2}}{2}}, \rho_{H_{1}}, \rho_{H_{2}}\in \ell ^{\frac{4}{3}}(\mathbb{Z})$.   Also, by symmetry 
\begin{equation*}
\mathbf{E} F_{N } ^ {(2)}\left( \frac{ \Vert DF_{N} ^{(1)} \Vert _{H} ^{2}}{2}-1\right)\to _{N\to \infty} \langle \rho_{\frac{H_{1}+H_{2}}{2}} \ast \rho_{H_{2}}, \rho_{\frac{H_{1}+H_{2}}{2}} \rangle _{\ell ^{2} (\mathbb{Z})}.
\end{equation*}
Next, denote by $ Y_{N} '$ the last component of the vector (\ref{25m-2}), i.e.
$$ Y'_{N}:=\sqrt{N} (\frac{1}{2}\langle DF_{N} ^{(1)}, DF_{N} ^{(2)}\rangle_{H}  - C_{1,2}).$$
By Lemma \ref{f12}
\begin{eqnarray*}
\mathbf{E} F _{N} ^{(1)}  Y'_{N} &\sim &\mathbf{E} F _{N} ^{(1)}  Y_{N}\sim C \frac{1}{N} \mathbf{E} I_{2} \left(\sum_{i=0} ^{N-1}(L_{k} ^{H_{1}} ) ^{\otimes 2}  \right) I_{2}\left( \sum_{j,k=0} ^{N-1} \rho_{\frac{H_{1}+H_{2}}{2}} (j-k) (\L ^{H_{1}}_{j} \otimes L ^{H_{2}}_{k} )\right) \\
&=& \frac{1}{N} \sum_{i,j,k=0} ^{N-1} \langle (L_{k} ^{H_{1}} ) ^{\otimes 2} , (L ^{H_{1}}_{j} \tilde{ \otimes } L ^{H_{2}}_{k} )\rangle  \\
&=& C \frac{1}{N} \sum _{i,j,k} \rho_{\frac{H_{1}+H_{2}}{2}} (j-k)\rho_{H_{1}} (i-j) \rho_{\frac{H_{1}+H_{2}}{2}} (i-k) \to _{N\to \infty } \langle \rho_{\frac{H_{1}+H_{2}}{2}} \ast \rho_{H_{1}}, \rho_{\frac{H_{1}+H_{2}}{2}} \rangle _{\ell ^{2} (\mathbb{Z})}
\end{eqnarray*}
and similarly
\begin{equation*}\mathbf{E} F_{N }^ {(2)}Y'_{N}\to _{N\to \infty } C\langle \rho_{\frac{H_{1}+H_{2}}{2}} \ast \rho_{H_{2}}, \rho_{\frac{H_{1}+H_{2}}{2}} \rangle _{\ell ^{2} (\mathbb{Z})}.
\end{equation*}
Finally, we can also prove that

\begin{eqnarray*}
\mathbf{E}\left( \frac{ \Vert DF_{N} ^{(1)}\Vert _{H} ^{2}}{2}-1\right)\left( \frac{ \Vert DF_{N} ^{(2)}\Vert _{H} ^{2}}{2}-1\right) &\sim&  C \frac{1}{N} \sum_{i,j,k,l} \rho_{H_{1}} (i-j) \rho_{H_{2}}(k-l) \langle (L_{i} ^{H_{1}}  \tilde{\otimes} L_{j} ^{H_{1}} ),( L_{k} ^{H_{2}}  \tilde{\otimes} L_{l} ^{H_{2}} ) \rangle \\
&=& C \frac{1}{N} \sum_{i,j,k.l}\rho_{H_{1}} (i-j) \rho_{H_{2}}(k-l)  \rho_{\frac{H_{1}+H_{2}}{2}} (k-i)\rho_{\frac{H_{1}+H_{2}}{2}} (l-j)\\
&& \to _{N \to \infty} \langle \rho_{H_{1}}\ast \rho_{\frac{H_{1}+H_{2}}{2}} , \rho_{H_{2}}\ast \rho_{\frac{H_{1}+H_{2}}{2}} \rangle _{\ell ^{2} (\mathbb{Z})}
\end{eqnarray*}
and 
\begin{equation*}
\mathbf{E}\left( \frac{ \Vert DF_{N} ^{(1)}\Vert _{H} ^{2}}{2}-1\right)Y'_{N}\sim \mathbf{E}\left( \frac{ \Vert DF_{N} ^{(1)}\Vert _{H}^{2}}{2}-1\right)Y_{N}\to _{N \to \infty}\langle \rho_{H_{1}}\ast \rho_{\frac{H_{1}+H_{2}}{2}} , \rho_{H_{1}}\ast \rho_{\frac{H_{1}+H_{2}}{2}} \rangle _{\ell ^{2} (\mathbb{Z})}
\end{equation*}

\begin{equation*}
\mathbf{E}\left( \frac{ \Vert DF_{N} ^{(2)}\Vert _{H} ^{2}}{2}-1\right)Y'_{N}\sim \mathbf{E}\left( \frac{ \Vert DF_{N} ^{(2)}\Vert _{H} ^{2}}{2}-1\right)Y_{N}\to _{N\to \infty}\langle \rho_{H_{2}}\ast \rho_{\frac{H_{1}+H_{2}}{2}} , \rho_{H_{2}}\ast \rho_{\frac{H_{1}+H_{2}}{2}} \rangle _{\ell ^{2} (\mathbb{Z})}.
\end{equation*}
So, condition (\ref{25m-4}) is fullfilled and the conclusion follows from Theorem \ref{m4m}.
\qed

\subsection{An example in a sum of finite Wiener chaoses}

Consider a sequence $I_{2} (f_{N})$ in a second Wiener chaos that satisfies condition {\bf  [C1\!\! ]}- {\bf  [C3\!\! ]} (take for example the sequence (\ref{vnn})). Consider a "perturbation" $G_{N}= I_{q} (g_{N}) $ with $q\geq 3$ and such that for $N\geq 1$
\begin{equation}
\label{25m-5}
\Vert g_{N} \Vert ^{2}_{H ^{\otimes q}} \leq C \frac{1}{N ^{1+\beta}}
\end{equation}
with some $\beta>0$. In particular, (\ref{25m-5}) implies that $G_{N}\to_{N\to \infty } 0$ in $L^{2} (\Omega)$. Define,
\begin{equation}
\label{25m-6}
F_{N} = I_{2}( f_{N})+ G_{N}, \hskip0.5cm N\geq 1.
\end{equation}
Clearly $F_{N}$ verifies {\bf  [C2\!\!  ]  }. Let us show that $F_{N}$ satisfies {\bf  [C1\!\! ] } with $\gamma_{N}= \sqrt{ Var \left( \frac{1}{2} \Vert D I_{2} (f_{N})\Vert _{H}^{2}\right)  }$ (which behaves as $C \frac{1}{\sqrt{N}}$ for $N$ large).  We have, by the product formula (\ref{prod})
\begin{eqnarray*}
\gamma_{N} ^{-1} \left( \langle DF_{N}, D(-L) ^{-1} F_{N} \rangle -1\right)
&=& \gamma_{N} ^{-1} \left( \frac{1}{2} \Vert D I_{2} (f_{N})\Vert _{H} ^{2}-1\right) + \frac{1}{q} \gamma_{N} ^{-1}\Vert DG_{N} \Vert _{H}^{2} \\
&&+ (2+q) \gamma_{N} ^{-1}\left( I_{q} (f_{N}\otimes _{1} g_{N})+ (q-1) I_{q-1} (f\otimes _{2} g) \right)
\end{eqnarray*}
Notice that, by (\ref{25m-5})
\begin{equation}
\label{25m-66}
\mathbf{E} \left( \gamma_{N}^{-1} \Vert DG_{N} \Vert _{H}^{2}\right)\leq CN ^{-\frac{1}{2}} \Vert g_{N} \Vert _{H ^{\otimes q}} ^{2} \leq C N ^{-\beta-\frac{1}{2}}\to_{N\to \infty }0.
\end{equation}
Also, since from (\ref{contra}) $\Vert f_{N}\otimes _{k} g_{N}\Vert _{H ^{\otimes (2+q-2k)}} \leq \Vert f_{N}\Vert _{H ^{\otimes 2}}\Vert g_{N}\Vert _{H^{\otimes q}}$ for  $k=1,2$ 
\begin{equation}
\label{25m-7}
\mathbf{E}\left( \gamma_{N}^{-1} I_{q} (f_{N}\otimes _{1} g_{N})\right)^{2} \leq CN \Vert f_{N} \otimes _{1} g\Vert_{H ^{\otimes q}} ^{2}\leq CN \Vert f_{N} \Vert ^{2} _{H ^{\otimes 2}}\Vert g_{N} \Vert ^{2} _{H ^{\otimes q}}\leq C N ^{-\beta}
\end{equation}
and
\begin{equation}
\label{25m-8}
\mathbf{E} \left( \gamma_{N}^{-1} I_{q-2} (f_{N}\otimes _{2} g_{N})\right)^{2} \leq CN \Vert f_{N} \otimes _{2} g\Vert_{H ^{\otimes (q-2)}} ^{2}\leq CN \Vert f_{N} \Vert ^{2} _{H^{\otimes 2}}\Vert g_{N} \Vert _{H ^{\otimes q}}^{2} \leq C N ^{-\beta}.
\end{equation}
 Relations (\ref{25m-66}), (\ref{25m-7}) and (\ref{25m-8}), together with the fact that $ \gamma_{N} ^{-1} \left( \frac{1}{2} \Vert D I_{2} (f_{N})-1\right) $ converges to a Gaussian law as $N\to \infty$, implies {\bf  [C1\!\! ]}. The same relations together with (\ref{hyper})  will give {\bf  [C3\!\! ]}.

\section{Appendix: Elements from Malliavin calculus}\label{app}

We briefly describe the tools from the analysis on Wiener space that we will need in our work. For complete presentations, we refer to \cite{N} or \cite{NPbook}. Consider $H$ a real separable Hilbert space and $(W(h), h \in H)$ an isonormal Gaussian process on a probability space $(\Omega, {\cal{A}}, P)$, which is a centered Gaussian family of random variables such that ${\bf E}\left[ W(\varphi) W(\psi) \right]  = \langle\varphi, \psi\rangle_{H}$. Denote by  $I_{n}$ the multiple stochastic integral with respect to
$B$ (see \cite{N}). This mapping $I_{n}$ is actually an isometry between the Hilbert space $H^{\odot n}$(symmetric tensor product) equipped with the scaled norm $\frac{1}{\sqrt{n!}}\Vert\cdot\Vert_{H^{\otimes n}}$ and the Wiener chaos of order $n$ which is defined as the closed linear span of the random variables $h_{n}(W(h))$ where $h \in H, \|h\|_{H}=1$ and $h_{n}$ is the Hermite polynomial of degree $n \in {\mathbb N}$
\begin{equation*}
h_{n}(x)=\frac{(-1)^{n}}{n!} \exp \left( \frac{x^{2}}{2} \right)
\frac{d^{n}}{dx^{n}}\left( \exp \left( -\frac{x^{2}}{2}\right)
\right), \hskip0.5cm x\in \mathbb{R}.
\end{equation*}
The isometry of multiple integrals can be written as follows: for $m,n$ positive integers,
\begin{eqnarray}
\mathbf{E}\left(I_{n}(f) I_{m}(g) \right) &=& n! \langle \tilde{f},\tilde{g}\rangle _{H^{\otimes n}}\quad \mbox{if } m=n,\nonumber \\
\mathbf{E}\left(I_{n}(f) I_{m}(g) \right) &= & 0\quad \mbox{if } m\not=n.\label{iso}
\end{eqnarray}
It also holds that
$I_{n}(f) = I_{n}\big( \tilde{f}\big)$ where $\tilde{f} $ denotes the symmetrization of $f$.
\\
We recall that any square integrable random variable which is measurable with respect to the $\sigma$-algebra generated by $W$ can be expanded into an orthogonal sum of multiple stochastic integrals
\begin{equation}
\label{sum1} F=\sum_{n=0}^\infty I_{n}(f_{n})
\end{equation}
where $f_{n}\in H^{\odot n}$ are (uniquely determined)
symmetric functions and $I_{0}(f_{0})=\mathbf{E}\left[  F\right]$.
\\\\
Let $L$ be the Ornstein-Uhlenbeck operator
\begin{equation*}
LF=-\sum_{n\geq 0} nI_{n}(f_{n})
\end{equation*}
if $F$ is given by (\ref{sum1}) and it is such that $\sum_{n=1}^{\infty} n^{2}n! \Vert f_{n} \Vert ^{2} _{{\cal{H}}^{\otimes n}}<\infty$.
\\\\
For $p>1$ and $\alpha \in \mathbb{R}$ we introduce the Sobolev-Watanabe space $\mathbb{D}^{\alpha ,p }$  as the closure of
the set of polynomial random variables with respect to the norm
\begin{equation*}
\Vert F\Vert _{\alpha , p} =\Vert (I -L) ^{\frac{\alpha }{2}} F \Vert_{L^{p} (\Omega )}
\end{equation*}
where $I$ represents the identity. We denote by $D$  the Malliavin  derivative operator that acts on smooth functions of the form $F=g(W(h_1), \dots , W(h_n))$ ($g$ is a smooth function with compact support and $h_i \in H$)
\begin{equation*}
DF=\sum_{i=1}^{n}\frac{\partial g}{\partial x_{i}}(W(h_1), \ldots , W(h_n)) h_{i}.
\end{equation*}
The operator $D$ is continuous from $\mathbb{D}^{\alpha , p} $ into $\mathbb{D} ^{\alpha -1, p} \left( H\right).$

We  recall  the product formula for multiple integrals.
It is well-known that for $f\in H^{\odot n}$ and $g\in H^{\odot m}$
\begin{equation}\label{prod}
I_n(f)I_m(g)= \sum _{r=0}^{n\wedge m} r! \left( \begin{array}{c} n\\r\end{array}\right) \left( \begin{array}{c} m\\r\end{array}\right) I_{m+n-2r}(f\otimes _r g)
\end{equation}
where $f\otimes _r g$ means the $r$-contraction of $f$ and $g$. This contraction is defined in the following way. Consider $(e_{j})_{j\geq 1}$ a complete orthonormal system in $H$ and let $f\in H^{\otimes n}$, $ g\in H^{\otimes m}$ be two symmetric functions with $n,m\geq 1$.  Then
\begin{equation}\label{f}
f= \sum_{j_{1}, .., j_{n} \geq 1} \lambda _{j_{1},.., j_{n}} e_{j_{1}} \otimes...\otimes e_{j_{n}}
\end{equation}
and
\begin{equation}\label{g}
g= \sum_{k_{1},.., k_{m}\geq 1} \beta _{ k_{1},.., k_{m}} e_{k_{1}}\otimes..\otimes e_{k_{m}} 
\end{equation}
where the coefficients $\lambda _{i}$ and $\beta _{j}$  satisfy $\lambda_{ j_{\sigma (1)},...j_{\sigma (n)}}= \lambda _{j_{1},.., j_{n}}$ and $ \beta _{ k_{\pi(1)},..., k_{\pi(m)}}= \beta _{ k_{1},.., k_{m}}$ for every permutation $\sigma $ of the set $\{1,..., n\}$ and for every permutation $\pi$ of the set $\{1,.., m\}$. Actually $\lambda _{j_{1},.., j_{n}}  = \langle f, e_{j_{1}} \otimes...\otimes e_{j_{n}}\rangle$ and $\beta _{ k_{1},.., k_{m}} =\langle g, e_{k_{1}}\otimes..\otimes e_{k_{m}} \rangle $ in (\ref{f}) and (\ref{g}).

If $f\in H^{\otimes n}$, $ g\in H^{\otimes m}$ are symmetric  given by (\ref{f}), (\ref{g}) respectively, then the contraction of order $r$ of $ f$ and $g$ is given by
\begin{eqnarray}
f\otimes _{r}g &=& \sum_{ i_{1},.., i_{r} \geq 1} \sum_{j_{1},..,j_{n-r}\geq 1} \sum_{k_{1},.., k_{m-r}\geq 1} \lambda _{i_{1},.., i_{r} , j_{1}, .., j_{n-r}}\beta _{i_{1},.., i_{r} , k_{1}, .., k_{m-r}}\nonumber\\
&&\times \left( e_{j_{1}}\otimes ..\otimes e_{j_{n-r}}\right) \otimes \left( e_{k_{1}}\otimes ..\otimes e_{k_{m-r}}\right) \label{contra}
\end{eqnarray}
for every $r=0,.., m\wedge n$. In particular $f\otimes _{0}g= f\otimes g$.
Note that $f\otimes _{r} g $ belongs to $H^{\otimes (m+n-2r)}$ for every $r=0,.., m\wedge n
$ and it is not in general symmetric. We will denote by $f\tilde{\otimes }_{r}g$ the symmetrization of $f\otimes _{r}g$.

Another important  property of  finite sums of multiple integrals is the hypercontractivity. Namely, if $F= \sum_{k=0} ^{n} I_{k}(f_{k}) $ with $f_{k}\in H ^{\otimes k}$ then
\begin{equation}
\label{hyper}
\mathbf{E}\vert F \vert ^{p} \leq C_{p} \left( \mathbf{E}F ^{2} \right) ^{\frac{p}{2}}.
\end{equation}
for every $p\geq 2$.

We denote by $D$  the Malliavin  derivative operator that acts on smooth functionals of the form $F=g(W(\varphi _{1}), \ldots , W(\varphi_{n}))$ (here $g$ is a smooth function with compact support and $\varphi_{i} \in H$ for $i=1,..,n$)
\begin{equation*}
DF=\sum_{i=1}^{n}\frac{\partial g}{\partial x_{i}}(W(\varphi _{1}), \ldots , W(\varphi_{n}))\varphi_{i}.
\end{equation*}
The operator $D$ can be extended to the space $\mathbb{D}^{1,2}$ since it is closable.

 The adjoint of $D$ is denoted by $\delta $ and is called the divergence (or
Skorohod) integral. Its domain ($Dom(\delta)$) coincides with the class of stochastic processes $u\in L^{2}(\Omega \times T)$ such that
\begin{equation*}
\left| \mathbf{E}\langle DF, u\rangle \right| \leq c\Vert F\Vert _{2}
\end{equation*}
for all $F\in \mathbb{D}^{1,2}$ and $\delta (u)$ is the element of $L^{2}(\Omega)$ characterized by the duality relationship
\begin{equation}\label{dua}
\mathbf{E}(F\delta (u))= \mathbf{E}\langle DF, u\rangle _{H}.
\end{equation}

The chain rule for the Malliavin derivative (see Proposition 1.2.4 in \cite{N}) will be used several times. If $\varphi: \mathbb{R}\to \mathbb{R}$ is a continuously differentiable function and $F\in \mathbb{D} ^ {1,2}$, then  $\varphi (F) \in \mathbb{D} ^ {1,2}$ and
\begin{equation}
\label{chain}
D\varphi(F)= \varphi ' (F) DF.
\end{equation}

\end{document}